\documentclass[10pt]{amsart}

\usepackage{graphicx}
\usepackage[top=26mm, bottom=26mm, left=26mm, right=26mm]{geometry}
\usepackage[dvips]{color}
\usepackage{bm}
\usepackage{marvosym}

\begin{document}

\begin{center}A level-set-based topology optimisation for acoustic-elastic coupled problems with a fast BEM-FEM solver \\
\vspace{10pt}
Hiroshi Isakari$^1$, Toyohiro Kondo$^2$, Toru Takahashi$^1$, and Toshiro Matsumoto$^1$ \\
1. Nagoya University, Japan \\
2. Toyota Industries Corporation, Japan 
\end{center}
\vspace{10pt}

\noindent
\textbf{Abstract:}  
This paper presents a structural optimisation method in three-dimensional acoustic-elastic coupled problems.
The proposed optimisation method finds an optimal allocation of elastic materials which reduces the sound level on some fixed observation points.
In the process of the optimisation, configuration of the elastic materials is expressed with a level set function, 
and the distribution of the level set function is iteratively updated with the help of the topological derivative.
The topological derivative is associated with state and adjoint variables which 
are the solutions of the acoustic-elastic coupled problems.
In this paper, the acoustic-elastic coupled problems are solved by a BEM-FEM coupled solver, 
in which the fast multipole method (FMM) and a multi-frontal solver for sparse matrices are efficiently combined.
Along with the detailed formulations for the topological derivative and the BEM-FEM coupled solver, 
we present some numerical examples of optimal designs of elastic sound scatterer to manipulate sound waves, 
from which we confirm the effectiveness of the present method.

\vspace{10pt}

\section{Introduction}
Computer simulations play an important role in modern product manufacturing process in various engineering fields.
The concept of computer aided engineering (CAE) is now widely accepted in some industries, 
which utilises numerical analysis to aid in tasks to evaluate the performance of engineering products.
With the help of CAE, the total cost and period for product developments have considerably been reduced.
In these days, use of computer simulations is not limited to performance evaluation, but is extended to design process.
As such an attempt, we can mention a structural optimisation, which is classified into sizing, shape and topology optimisation \cite{bendsoe1988generating}.
The topology optimisation is considered as the most powerful design method in structural optimisations since it can design not only the shape but also the topology of devices, 
i.e., the topology optimisation allows an nucleation of a new material and/or hole in its process. 
Hence, the obtained optimal design by the topology optimisation is less affected by an initial guess than the other structural optimisations.

After a pioneering work in \cite{bendsoe1988generating}, 
the topology optimisation is intensively studied mainly in the field of structural mechanics in order to design light but stiff structural members~\cite{bendsoe1988generating,suzuki1991homogenization,yamada2010topology}.
Recently, one of the main interests in the topology optimisation community is to widen the applicability of the topology optimisation to problems in various engineering fields other than structural mechanics, 
such as thermal problems~\cite{yamada2011level,jing2014topological,li2016generating}, fluid problems~\cite{yaji2014topology,papoutsis2015continuous,yaji2016topology}, elastodynamic problems~\cite{jensen2007toplogy,sigmund2003systematic,ma1995topological}, electromagnetic wave problems~\cite{yoo2000topology,abe2013level_en,kourogi2013electromagnetic_en,isakari2016multi,heo2016dielectric}, and so on.
There are also many efforts to extend the topology optimisation into various design problems in acoustics such as an acoustic horn to maximise sound level~\cite{wadbro2006topology}, 
an acoustic metamaterial which realises a material with negative effective bulk modulus~\cite{lu2013topology}, 
and a poroelastic sound-proofing material~\cite{lee2007optimal,yamamoto2009topology}.
We think, however, the applicability of these existing methods for industrial designs is still limited 
because these topology optimisations use the finite element method (FEM) to solve boundary value problems involved in sensitivity analysis.
Since the acoustic problem is often defined in an unbounded domain, the unbounded domain is approximated with a large one in FEM, which leads to an unexpected large scale problem.
Further, an artificial boundary condition such as perfect matched layer (PML) is required to make sure that the scattered wave does not reflect on the truncated boundary.

On the other hand, when the boundary element method (BEM) is employed to solve wave scattering problems, 
only the boundary of the domain is needed to be discretised, 
which considerably reduces the number of elements. 
Also, the scattered fields by the BEM automatically satisfy the radiation condition, i.e., 
no artificial boundary condition is required to deal with the unbounded domain with the BEM.
Thus, the BEM is more suitable for topology optimisations in wave problems than the FEM.
As pioneering works on BEM-based topology optimisations, we can mention Abe et al~\cite{abe2007boundary} and Du and Olhoff~\cite{du2007minimization}, 
In the first one, they have solved a two-dimensional shape optimisation problem by the BEM to design a sound barrier. 
In the second one, they have solved a topology optimisation problem for a noise reduction device from vibrating structures, 
in which they use an approximated boundary integral formulation for high frequency problems.
The applicability of these BEM-based method is limited to either two-dimensional problem~\cite{abe2007boundary} or high frequency problem~\cite{du2007minimization} because a naive BEM for three-dimensional realistic scale problems is too expensive.
It is inevitable to accelerate the BEM by, for example, fast multipole method~(FMM)~\cite{rokhlin1985rapid,greengard1987fast}, $\mathcal H$ matrix algebra \cite{bebendorf2008hierarchical} and fast direct solver~\cite{martinsson2005fast} 
for topology optimisation in sound problems.

In order to realise a topology optimisation for three-dimensional realistic design problems of wave devices,
we have been investigating level-set-based topology optimisations with the BEM accelerated by the FMM.
In our methodology, a candidate for optimal configuration is expressed with a level set function 
which is iteratively updated with the help of the topological derivative~\cite{isakari2014topology,bonnet2007fm,novotny2003topological,carpio2008solving} to find an optimal distribution of sound scatterers.
In \cite{isakari2014topology}, we have investigated a topology optimisation for rigid materials to minimise sound pressure on some observation points.
We have extended the methodology to find an optimal allocation of sound absorbers in \cite{kondo2014acoustic_en}, 
in which a sound absorbing material is modelled with the impedance boundary condition.
The impedance boundary condition is, however, not appropriate to model sound absorbing material in some applications.
For example, in analysis with the impedance boundary condition, penetrated sound waves in the sound absorbing materials cannot explicitly be observed, and vibrations in the sound absorber itself are neglected.

In this study, to further enhance the applicability of our methodology, we present a level-set-based topology optimisation in acoustic-elastic coupled problems, 
with which the vibrations of sound scatterers made of elastic materials are explicitly evaluated.
In order to solve the acoustic-elastic coupled problem, we adopt a BEM-FEM solver which solves the acoustic and elastic field by BEM and FEM, respectively.
This choice is reasonable since, with our settings, the elastic material is in a bounded domain while the acoustic host matrix is unbounded.
Although the acoustic-elastic coupled problem can appropriately be solved by the BEM~\cite{isakari2012periodic_gakui,kimeswenger2014coupled}, 
we use the BEM-FEM solver~\cite{everstine1990coupled} since the solver can naturally be extended to deal with elastic material other than the isotropic one such as anisotropic material and Biot's poroelastic material~\cite{biot1956theory,biot1956theory2}.
So far, some fast techniques~\cite{fischer2005fast} for the BEM-FEM solver are proposed.
We here propose another acceleration technique for the BEM-FEM coupled solver, 
in which FMM and a multi-frontal solver for sparse matrices are efficiently combined.

The rest of the paper is organised as follows.
After presenting the statement of the acoustic-elastic coupled problem and related optimisation problem in Section 2.1, 
we derive the relevant topological derivatives in Section 2.2.
In Section 2.3, we propose a new fast method with the FMM and a multi-frontal solver to solve algebraic equations stemmed from the BEM-FEM coupling for the acoustic-elastic coupled problem.
After briefly reviewing a formulation and algorithm of the present optimisation method in Sections 2.4 and 2.5, 
we present some numerical examples which verify the efficiency of the proposed methods in Section 3.
Specifically, we check the computational cost for the present BEM-FEM coupled solver in Section 3.1, 
numerically verify the topological derivatives in Section 3.2, 
and present two optimal designs of elastic sound scatterers in Section 3.3.
In Section 4, we conclude the paper, and discuss the remaining issues to be addressed in the future. 

\section{Formulations}
\subsection{Statement of the optimisation problem in acoustic-elastic coupled problems}
We consider a three-dimensional acoustic-elastic coupled problem in which
an incident sound wave from sources on $\boldsymbol{x}^\mathsf{src}_i~(i=1,...,M^\mathsf{src})$, where $M^\mathsf{src}$ is the number of the sources, is scattered by elastic scatterers filled in a bounded domain $\Omega^\mathsf{c}$.
The sound field in $\Omega:=\mathbb{R}^3\setminus\overline{\Omega^\mathsf{c}}$ and the transmitted elastic field in $\Omega^\mathsf{c}$ are governed by the following boundary value problem:
\begin{align}
 p_{,jj}\left(\boldsymbol{x}\right)+k^2_\mathsf{f} p\left(\boldsymbol{x}\right)+\sum_{m=1}^{M^\mathsf{src}} A^{\mathsf{src}}_m \delta\left(\boldsymbol{x}-\boldsymbol{x}_{m}^{\mathsf{src}}\right)=0& \ \ \ \ \ \ \ \ \boldsymbol{x}\in\Omega,   \label{hel}\\
\sigma_{ji,j} \left(\boldsymbol{x}\right)+\rho_\mathsf{s} \omega^2 u_i \left(\boldsymbol{x}\right)=0& \ \ \ \ \ \ \ \ \boldsymbol{x}\in\Omega^\mathsf{c},   \label{ela}\\
t_i(\boldsymbol{x})+ p\left(\boldsymbol{x}\right)n_i\left(\boldsymbol{x}\right)=0& \ \ \ \ \ \ \ \ \boldsymbol{x}\in\Gamma:=\overline{\partial\Omega\cap\partial\Omega^\mathsf{c}},  \label{bc_power}\\
q\left(\boldsymbol{x}\right):=\frac{\partial p(\boldsymbol{x})}{\partial n(\boldsymbol{x})}=\rho_\mathsf{f}\omega^2 u_i\left(\boldsymbol{x}\right)n_i\left(\boldsymbol{x}\right)&  \, \ \ \ \ \ \ \  \boldsymbol{x}\in\Gamma, \label{bc_velo}\\
\left|\boldsymbol{x}\right| \left( q\left(\boldsymbol{x}\right)- \mathsf{i}k_\mathsf{f} p\left(\boldsymbol{x}\right)\right)\rightarrow 0 & \ \ \ \ \ \ \ \ \mathrm{as} \  \ |\boldsymbol{x}|\rightarrow \infty,\label{rad}
\end{align}
where $A^\mathsf{src}_m$ is the intensity of $m$-th sound source, and $\delta$ is the Dirac delta,
$p$ is the sound pressure, $\sigma_{ij}$ is the stress, $u_i$ is the displacement, and $t_i$ is the traction defined as
\begin{align}
t_i(\boldsymbol{x}) =C_{ijk\ell} u_{k,\ell j}(\boldsymbol{x})n_j(\boldsymbol{x}), \label{traction}
\end{align}
where $\boldsymbol{n}(\boldsymbol{x})$ is the exterior normal vector with respect to $\Omega^\mathsf{c}$ on $\boldsymbol{x}\in\Gamma$, 
and $C_{ijk\ell}$ is the elastic tensor which, for the case that an isotropic material is concerned, has the following representation:
\begin{align}
C_{ijk\ell}=\lambda\delta_{ij}\delta_{k\ell}+\mu\delta_{ik}\delta_{j\ell}+\mu\delta_{i\ell}\delta_{jk}, \label{cijkl}
\end{align}
where $\lambda$ and $\mu$ are the Lam\'e constants, and $\delta_{ij}$ is the Kronecker delta.
Also, $\rho_\mathsf{f}$ and $\rho_\mathsf{s}$ are the densities for acoustic and elastic materials, respectively, and $k_\mathsf{f}$ is the wave number for the acoustic wave defined as
\begin{align}
 k_\mathsf{f}=\omega\sqrt{\frac{\rho_\mathsf{f}}{\Lambda_\mathsf{f}}},
 \label{kf}
\end{align}
where $\Lambda_\mathsf{f}$ is the bulk modulus of the acoustic material, 
and $\omega$ is the frequency with which the time dependency of physical quantities is assumed to be $e^{-\mathsf{i}\omega t}$.

The boundary value problem in Eqs.~\eqref{hel}--\eqref{rad} is not uniquely solvable for certain frequencies, called eigenfrequencies~\cite{kimeswenger2014coupled,luke1995fluid}.
With the eigenfrequency $\omega$, there exists an non-trivial solutions of the homogeneous boundary value problem, i.e. for $A_m^\mathsf{src}=0~(m=1,...,M^\mathsf{src})$ in \eqref{hel}.
For a real eigenfrequency called the ``Jones frequency'', the corresponding non-trivial solution satisfies $p=0$ in $\Omega$ and $u_in_i=0$, $t_i=0$ on $\Gamma$.
It is also known that complex eigenfrequencies with negative imaginary part may exist,
which may affect the accuracy of numerical methods when the relevant frequency is close (even when not identical) to one of the eigenfrequencies.
We henceforth assume that the frequency $\omega$ is real and far away from any of the eigenfrequencies.

\begin{figure}[h]
 \begin{center}
  \includegraphics[scale=0.35]{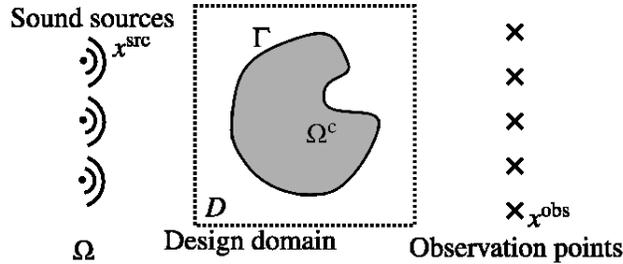}
  \caption{Settings for the topology optimisation in acoustic-elastic coupled problem.}
  \label{01formulation_optimisation}
 \end{center}
\end{figure}

Our optimisation problem is defined as to find an optimal distribution of elastic material(s) $\Omega^\mathsf{c}\subset D$ (Figure~\ref{01formulation_optimisation})
which minimises the following objective function $J$ defined with a functional $f$:
\begin{align}
 J=\sum_{m=1}^{M^\mathsf{obs}} f\left(p\left(\boldsymbol{x}_m^\mathsf{obs}\right)\right), \label{obj}
\end{align}
where $D$, which is so called design domain, is bounded, 
and $\boldsymbol{x}_m^\mathsf{obs}\notin D$ is $m$-th observation point on which the sound level is evaluated, 
and $M^\mathsf{obs}$ is the number of the observation points.

\subsection{The topological derivative for acoustic-elastic coupled problems}\label{sec:td}
In this subsection, we present the topological derivative for the objective function in Eq.~\eqref{obj}.
We here derive the topological derivative ${\mathcal T}_\Omega$ which characterises the sensitivity of $J$ to an appearance of an infinitesimal spherical elastic material $\Omega_\varepsilon$ in the acoustic matrix $\Omega$.
The topological derivative ${\mathcal T}_{\Omega^\mathsf{c}}$ with respect to appearance of an acoustic material $\Omega_\varepsilon$ in the elastic inclusion $\Omega^\mathsf{c}$ can similarly be obtained.

\begin{figure}[h]
 \begin{center}
  \includegraphics[scale=0.35]{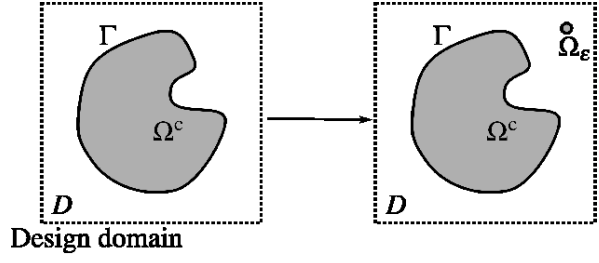}
  \caption{An infinitesimal elastic inclusion $\Omega_\varepsilon$ is introduced in the design domain $D$.}
  \label{02topology}
 \end{center}
\end{figure}

Let us assume that an infinitesimal spherical elastic material $\Omega_\varepsilon$, whose elastic properties are identical to those of $\Omega^\mathsf{c}$, appears in $D$ (Figure~\ref{02topology}). 
We henceforth denote the centre and the radius of the infinitesimal sphere $\Omega_\varepsilon$ as $\boldsymbol{x}^0$ and $\varepsilon$, respectively.
Due to the appearance of the small elastic inclusion $\Omega_\varepsilon$, 
the functions $p$, $q$, $u_i$, $t_i$, $\sigma_{ij}$ in Eqs.~\eqref{hel}--\eqref{rad} suffer from perturbations which are henceforth denoted as in Table~\ref{table1}.
\begin{table}[h]
 \begin{center}
  \caption{Perturbed physical quantities due to the appearance of the elastic inclusion $\Omega_\varepsilon$.}
  \label{table1}
  \begin{tabular}{ccc}
   \hline
   in $\Omega\setminus\overline{\Omega_\varepsilon}$ & in $\Omega_2$ & in $\Omega_\varepsilon$ \\
   \hline
   sound pressure: $p+\delta p$ & displacement $u_i+\delta u_i$ & displacement $\hat{u}_i$ \\
   sound flux: $q+\delta q$ & stress $\sigma_{ij}+\delta\sigma_{ij}$ & stress $\hat{\sigma}_{ij}$ \\
   -- & traction: $t_i+\delta t_i$ & traction: $\hat{t}_i$\\
   \hline
  \end{tabular}
 \end{center}
\end{table}
The perturbations are governed by the following boundary value problem:
\begin{align}
\delta p_{,jj}\left(\boldsymbol{x}\right) + k^2_\mathsf{f} \delta p\left(\boldsymbol{x}\right) = 0&  \ \ \ \ \ \boldsymbol{x}\in\Omega \setminus \overline{\Omega}_\varepsilon, \label{danseidelhel} \\
\delta \sigma_{ji,j} \left(\boldsymbol{x}\right)+\rho_{\mathsf{s}} \omega^2 \delta u_i \left(\boldsymbol{x}\right)=0& \ \ \ \ \ \boldsymbol{x}\in\Omega^\mathsf{c}, \label{delela}\\
\delta t_i\left(\boldsymbol{x}\right)+ \delta p\left(\boldsymbol{x}\right)n_i(\boldsymbol{x})=0 & \ \ \ \ \ \boldsymbol{x}\in\Gamma, \label{delbc_power}\\
\delta q\left(\boldsymbol{x}\right)=\rho_\mathsf{f} \omega^2 \delta u_i\left(\boldsymbol{x}\right)n_i(\boldsymbol{x}) & \ \ \ \ \ \boldsymbol{x}\in\Gamma, \label{delbc_velo}\\
\delta \hat{\sigma}_{ji,j} \left(\boldsymbol{x}\right)+\rho_{\mathsf{s}} \omega^2 \delta \hat{u}_i \left(\boldsymbol{x}\right)=0& \ \ \ \ \ \boldsymbol{x}\in\Omega_\varepsilon, \label{delelasig}\\
\hat{t}_i\left(\boldsymbol{x}\right)+ \left(p+\delta p\right)\left(\boldsymbol{x}\right)n_i(\boldsymbol{x})=0 & \ \ \ \ \ \boldsymbol{x}\in\Gamma_\varepsilon, \label{delbc_power2}\\
 \left(q+\delta q\right)\left(\boldsymbol{x}\right)=\rho_\mathsf{f} \omega^2 \hat{u}_i\left(\boldsymbol{x}\right)n_i(\boldsymbol{x}) & \ \ \ \ \ \boldsymbol{x}\in\Gamma_\varepsilon, \label{delbc_velo2}\\
\left|\boldsymbol{x}\right| \left( \delta q\left(\boldsymbol{x}\right)- \mathsf{i}k_\mathsf{f} \delta p\left(\boldsymbol{x}\right)\right)\rightarrow 0& \ \ \ \ \ \mathrm{as} \ |\boldsymbol{x}|\rightarrow \infty, \label{danseidelsommarfeld}
\end{align}
where the exterior normal vector on $\Gamma_\varepsilon:=\overline{\partial\Omega_\varepsilon\cap\partial\Omega}$ is defined with respect to $\Omega_{\varepsilon}$.
The objective function in \eqref{obj} also suffers from a perturbation $\delta J$ due to the appearance of $\Omega_\varepsilon$ as 
\begin{align}
\delta J={\Re}\left[
 \sum_{m=1}^{M^\mathsf{obs}}\frac{\partial f\left(p\left(\boldsymbol{x}_m^{\mathsf{obs}}\right)\right)}{\partial p} \delta  p\left(\boldsymbol{x}_m^{\mathsf{obs}}\right)\right].
\label{danseideltaJ}
\end{align}
Note that the direct evaluation of $\delta J$ with Eq.~\eqref{danseideltaJ} is impractical 
since it involves the perturbation of the sound pressure $\delta p$ on all observation points $\boldsymbol{x}^\mathsf{obs}_m$ 
which is the solution of the boundary value problem in Eqs.~\eqref{danseidelhel}--\eqref{danseidelsommarfeld}.

In this paper, we use the adjoint variable method to evaluate $\delta J$ without going through $\delta p(\boldsymbol{x}_m^\mathsf{obs})$.
The adjoint problem is defined as follows:
\begin{align}
\tilde{p}_{,jj}\left(\boldsymbol{x}\right) + k^2_\mathsf{f} \tilde{p}\left(\boldsymbol{x}\right)+\sum_{m=1}^{M^\mathsf{obs}}\frac{\partial f\left(p\left(\boldsymbol{x}_m^{\mathsf{obs}}\right)\right)}{\partial p} \delta\left(\boldsymbol{x}-\boldsymbol{x}_m^{\mathsf{obs}}\right) = 0&  \ \ \ \ \ \boldsymbol{x}\in\Omega, \label{danseiajohel} \\
\tilde{\sigma}_{ji,j}\left(\boldsymbol{x}\right)+\rho_\mathsf{s} \omega^2 \tilde{u}_i\left(\boldsymbol{x}\right)=0&  \ \ \ \ \ \boldsymbol{x}\in\Omega^\mathsf{c}, \label{danseiajoela} \\
\tilde{t}_i\left(\boldsymbol{x}\right)+ \tilde{p}\left(\boldsymbol{x}\right)n_i(\boldsymbol{x})=0 & \ \ \ \ \ \boldsymbol{x}\in\Gamma, \label{danseiajobc_power}\\
\tilde{q}\left(\boldsymbol{x}\right):=\frac{\partial\tilde{p}(\boldsymbol{x})}{\partial n}=\rho_\mathsf{f} \omega^2 \tilde{u}_i\left(\boldsymbol{x}\right)n_i(\boldsymbol{x}) & \ \ \ \ \ \boldsymbol{x}\in\Gamma, \label{danseiajobc_velo}\\
|\boldsymbol{x}|\left( \tilde{q}(\boldsymbol{x})- \mathsf{i}k_\mathsf{f} \tilde{p}(\boldsymbol{x})\right)\rightarrow0& \ \ \ \ \ \mathrm{as} \ |\boldsymbol{x}|\rightarrow\infty, \label{danseiajosommarfeld}
\end{align} 
where $\tilde{p}$, $\tilde{q}$, $\tilde{u}_i$, $\tilde{\sigma}_{ij}$ and $\tilde{t}_i$ are the adjoint sound pressure, the adjoint sound flux, the adjoint displacement, the adjoint stress and the adjoint traction, respectively.
According to the reciprocity of the state variable $p$ and the adjoint variable $\tilde{p}$ in $\Omega_\varepsilon$, we have the following identity:
\begin{align}
 \int_{\Gamma_\varepsilon} \left(\tilde{p}q -p\tilde{q}\right) \mathrm{d}\Gamma=0. \label{29_danseisouhan2}
\end{align}
A similar procedure to $\hat{u}_i$ and $\tilde{p}$ in $\Omega_\varepsilon$ together with Eqs.\eqref{delbc_power2}, \eqref{delbc_velo2}, \eqref{danseiajohel} and \eqref{29_danseisouhan2} gives the following identity:
\begin{align}
 \int_{\Gamma_\varepsilon}\left(\tilde{q}\delta p - \tilde{p} \delta q\right)\mathrm{d}\Gamma=&
\int_{\Omega_\varepsilon}\left(
 \Lambda_\mathsf{f}\hat{u}_{i,i}\tilde{p}_{,jj}
-\tilde{p}_{,ij}\hat{\sigma}_{ji}
+\rho_\mathsf{s} \omega^2 \tilde{p}_{,i}\hat{u}_i
-\rho_\mathsf{f}\omega^2 \hat{u}_i\tilde{p}_{,i}
\right) \mathrm{d}\Omega.
\label{32_danseisouhan3-2}
\end{align}
We also have the following reciprocal relation between $\delta u_i$ and $\tilde{p}$ in $\Omega^\text{c}$ combined with the boundary condition in Eq.~\eqref{delbc_power}:
\begin{align}
 \int_\Gamma \left(\delta u_i \tilde{t}_i+\delta p \tilde{u}_i n_i\right) \mathrm{d}\Gamma =0. \label{28_danseisouhan1-1}
\end{align}
With the reciprocal theorem between $p$ and $\tilde{p}$ in $\Omega\setminus\overline{\Omega_\varepsilon}$, 
and Eqs. \eqref{delbc_velo}, \eqref{danseiajobc_power}, \eqref{danseiajobc_velo}, \eqref{32_danseisouhan3-2} and \eqref{28_danseisouhan1-1}, 
we can evaluate $\delta J$ as follows:
\begin{align}
 \delta J
 =&{\Re}\left[\int_{\Omega_\varepsilon}\left(
 \Lambda_\mathsf{f}\hat{u}_{i,i}\tilde{p}_{,jj}
 -\tilde{p}_{,ij}\hat{\sigma}_{ji}
 +\rho_\mathsf{s} \omega^2 \tilde{p}_{,i}\hat{u}_i
 -\rho_\mathsf{f}\omega^2 \hat{u}_i\tilde{p}_{,i}
 \right) \mathrm{d}\Omega \right].
\label{35_danseisouhan4-2}
\end{align}
Note that the expression in Eq.~\eqref{35_danseisouhan4-2} does not involve the perturbations of the state variables on the observation points.
The expression \eqref{35_danseisouhan4-2} can further be simplified with the help of the Gauss theorem as
\begin{align}
\delta J
=&{\Re}\left[-\int_{\Gamma_\varepsilon} \left(\rho_\mathsf{f}\omega^2 \hat{u}_r \tilde{p}+\hat{t}_i \tilde{p}_{,i}\right) \mathrm{d}\Gamma \right],
\label{danseisouhan4-5}
\end{align}
where $\hat{u}_r:=\hat{u}_in_i$ is the radial component of the displacement on $\Gamma_\varepsilon$.

In the following, we evaluate the asymptotic behaviour of $\delta J$ in Eq.~\eqref{danseisouhan4-5} as $\varepsilon\rightarrow 0$.
To this end, $\tilde{p}$ and its gradient are respectively expanded as 
\begin{align}
 \tilde{p}(\boldsymbol{x})=\tilde{p}^0 +\varepsilon \tilde{p}^{0}_{,j}n_j(\boldsymbol{x})+o\left(\varepsilon\right) & \ \ \ \ \ \boldsymbol{x}\in\Gamma_\varepsilon, \label{kyutentilp}\\
 \tilde{p}_{,i}(\boldsymbol{x})=\tilde{p}^0 _{,i} +\varepsilon \tilde{p}^{0}_{,ij}n_j(\boldsymbol{x})+o\left(\varepsilon\right) & \ \ \ \ \ \boldsymbol{x}\in\Gamma_\varepsilon, \label{kyutentilpi}
\end{align}
where $f^0$ denotes $f^0=f(\boldsymbol{x}^0)$.
The state variables $p$ and $q$ can also be expanded as
\begin{align}
 p(\boldsymbol{x})=p^0+\varepsilon p^0_in_i(\boldsymbol{x})+o\left(\varepsilon\right) & \ \ \ \ \ \boldsymbol{x}\in\Gamma_\varepsilon, \label{kyutenp} \\
 q(\boldsymbol{x})=p^0_{,i} n_i(\boldsymbol{x})+\varepsilon p^0_{,ij}n_i(\boldsymbol{x})n_j(\boldsymbol{x})+o\left(\varepsilon\right) & \ \ \ \ \ \boldsymbol{x}\in\Gamma_\varepsilon, \label{kyutenq}
\end{align}
respectively.
The asymptotic expansions of the radial displacement $\hat{u}_r$ and the traction $\hat{t}_i$ on $\Omega_\varepsilon$ in \eqref{danseisouhan4-5} are then evaluated.
Although $\hat{u}_r$ and $\hat{t}_i$ are the solution of the boundary value problem in \eqref{danseidelhel}--\eqref{danseidelsommarfeld}, 
it is sufficient to solve the following approximated one~\cite{bonnet2007fm}:
\begin{align}
\delta p_{,jj}\left(\boldsymbol{x}\right) + k^2_\mathsf{f} \delta p\left(\boldsymbol{x}\right) = 0&  \ \ \ \ \boldsymbol{x} \in \Omega^\mathsf{c}_\varepsilon:=\mathbb{R}^3\setminus\overline{\Omega^\mathsf{c}}, \label{danzehel} \\
\hat{\sigma}_{j,ij} \left(\boldsymbol{x}\right)+\rho_{\mathsf{s}} \omega^2 \hat{u}_i \left(\boldsymbol{x}\right)=0& \ \ \ \ \ \boldsymbol{x} \in \Omega_\varepsilon,  \label{danzenela}\\
\left(q+\delta q\right)\left(\boldsymbol{x}\right)=\rho_\mathsf{f} \omega^2 \hat{u}_i\left(\boldsymbol{x}\right)n_i(\boldsymbol{x}) & \ \ \ \ \ \boldsymbol{x} \in \Gamma_\varepsilon, \label{danzenbc_velo}\\
\delta t_i\left(\boldsymbol{x}\right)+ \left(p+\delta p\right)\left(\boldsymbol{x}\right)n_i(\boldsymbol{x})=0& \ \ \ \ \ \boldsymbol{x}\in\Gamma_\varepsilon, \label{danzenbc_power}\\
|\boldsymbol{x}| \left( \delta q(\boldsymbol{x})- \mathsf{i}k_\mathsf{f} \delta p(\boldsymbol{x}) \right)\rightarrow 0 & \ \ \ \ \ \mathrm{as} \ |\boldsymbol{x}| \rightarrow \infty, 
\label{danzensommarfeld}
\end{align} 
since the asymptotic behaviour as $\varepsilon\rightarrow 0$ is now concerned.
For the sake of reference, we rewrite the boundary conditions \eqref{danzenbc_velo} and \eqref{danzenbc_power} as follows:
\begin{align}
 q\left(\boldsymbol{x}\right)+\delta q\left(\boldsymbol{x}\right)
=\rho_\mathsf{f} \omega^2 \hat{u}_r\left(\boldsymbol{x}\right)  \label{kyuten1}  \,\ \ \ \ \ \ \ \ \boldsymbol{x}\in\Gamma_\varepsilon, \\ 
 p\left(\boldsymbol{x}\right)+\delta p\left(\boldsymbol{x}\right)+\hat{\sigma}_{rr}\left(\boldsymbol{x}\right)=0  \ \ \ \ \ \ \ \ \ \ \ \ \boldsymbol{x}\in\Gamma_\varepsilon,  \label{kyuten2} \\
 \hat{\sigma}_{r\theta}\left(\boldsymbol{x}\right)=\hat{\sigma}_{r\phi}\left(\boldsymbol{x}\right)=0 \ \ \, \ \ \ \ \ \ \ \ \ \ \ \ \ \ \ \ \ \ \boldsymbol{x}\in\Gamma_\varepsilon, \label{kyuten3}
\end{align}
where $\hat{\sigma}_{rr}$, $\hat{\sigma}_{r\theta}$ and $\hat{\sigma}_{r\phi}$ denote the polar representations of the stress $\hat{\sigma}_{ij}$.
The solutions of the boundary value problem \eqref{danzehel}, \eqref{danzenela}, \eqref{danzensommarfeld}, \eqref{kyuten1}, \eqref{kyuten2} and \eqref{kyuten3} can be written in terms of spherical functions \cite{eringen1978elastodynamics} as follows:
\begin{align}
 \delta p\left(\boldsymbol{x}\right)&=\sum_{n=0}^{\infty}\sum_{m=-n}^{n}d_n^m {h_n^{\left(1\right)}} \left(k_\mathsf{f}r\right)P_n^m(\cos\theta)e^{\mathsf{i}m\phi} \ \ \ \ \ \boldsymbol{x}\in\Omega_\varepsilon^\mathsf{c},
\label{kyutendelp} \\
 \delta q\left(\boldsymbol{x}\right)&=\sum_{n=0}^{\infty}\sum_{m=-n}^{n}d_n^m \frac{\partial h_n^{\left(1\right)}\left(k_\mathsf{f}r\right)}{\partial r}P_n^m(\cos\theta)e^{\mathsf{i}m\phi} \ \ \ \ \ \boldsymbol{x}\in\Omega_\varepsilon^\mathsf{c},
\label{kyutendelq}\\
 \hat{u}_{r}(\boldsymbol{x})&=\sum_{n=0}^{\infty}\sum_{m=-n}^{n}
 \frac{1}{r}\left(a_{n}^{m} U_{1}^{n}(r)-\frac{c_n^m}{k_\mathsf{T}}U_{3}^{n}(r)\right)P_n^m(\cos\theta)e^{\mathsf{i}m\phi} \ \ \ \ \ \boldsymbol{x}\in\Omega_\varepsilon,\label{kyutenhatur}  \\
 \hat{\sigma}_{rr}(\boldsymbol{x})&=\sum_{n=0}^{\infty}\sum_{m=-n}^{n}\frac{2\mu}{r^2}\left(a_n^mT_{11}^n(r)-\frac{c_n^m}{k_\mathsf{T}}T_{13}^n(r)\right)P_n^m(\cos\theta)e^{\mathsf{i}m\phi} \ \ \ \ \ \boldsymbol{x}\in\Omega_\varepsilon,\label{kyutenhatsigmaij}\\
 \hat{\sigma}_{r\theta}\left(\boldsymbol{x}\right)&=
 \sum_{n=0}^{\infty}\sum_{m=-n}^{n}\frac{2\mu}{r^2}
\left[
  a_n^mT_{41}^n(r) \left(n\cot{\theta}P_n^m\left(\cos{\theta}\right)-\frac{n+m}{\sin{\theta}}P_{n-1}^m\left(\cos{\theta}\right)\right)
 \right. \nonumber \\
 &+b_n^mT_{42}(r)\frac{\mathsf{i}m}{\sin{\theta}}P_n^m\left(\cos{\theta}\right) \nonumber \\
 &\left.+c_n^mT_{43}^n(r) \left(n\cot{\theta}P_n^m\left(\cos{\theta}\right)-\frac{n+m}{\sin{\theta}}P_{n-1}^m\left(\cos{\theta}\right)\right)\right]e^{\mathsf{i}m\phi} \ \ \ \ \ x\in\Omega_\varepsilon,\label{hatsigmartheta} \\
 \hat{\sigma}_{r\phi}\left(\boldsymbol{x}\right)&=
 \sum_{n=0}^{\infty}\sum_{m=-n}^{n}\frac{2\mu}{r^2}
\left(
  a_n^mT_{51}^n(r) \frac{\mathsf{i}m}{\sin{\theta}}P_n^m\left(\cos{\theta}\right) \right. \nonumber \\
 &\left.-b_n^mT^n_{52}(r)\frac{\mathsf{i}m}{\sin{\theta}}\left(n\cos{\theta}P_n^m\left(\cos{\theta}\right)-\left(n+m\right)P_{n-1}^m\left(\cos{\theta}\right)\right) \right. \nonumber \\
 &\left.+\frac{c_n^m}{k_\mathsf{T}}T_{53}^n(r) \frac{\mathsf{i}m}{\sin{\theta}}P_n^m\left(\cos{\theta}\right)
\right)e^{\mathsf{i}m\phi} \ \ \ \ \ x\in\Omega_\varepsilon,
\label{hatsigmarphi}
\end{align}
where $h_n^{(1)}$ is the $n$-th spherical Hankel function of the first kind, $P_n^m$ is the associated Legendre function defined as follows:
\begin{align}
P_n^m(x)=(1-x^2)^{m/2}\frac{d^m}{dx^m}P_n(x)\hspace{10pt}(m\ge 0),\\
P_n^{-m}(x)=(-1)^m\frac{(n-m)!}{(n+m)!}P_n^m(x)\hspace{10pt}(m\ge 0),
\end{align}
where $P_n$ is the Legendre polynomial, and $(r,\theta,\phi)$ represents the spherical coordinate of the point $\boldsymbol{x}$.
Also, $U_1^n$, $U_3^n$, $T_{11}^n$, $T_{13}^n$, $T_{41}^n$, $T_{42}^n$, $T_{43}^n$, $T_{51}^n$, $T_{52}^n$ and $T_{53}^n$ are the functions defined as follows:
\begin{align}
 U_{1}^n(r)=&nj_n\left(k_\mathsf{L}r\right)-k_\mathsf{L}r j_{n+1}\left(k_\mathsf{L}r\right), \\
 U_{3}^n(r)=&n\left(n+1\right)j_n\left(k_\mathsf{T}r\right), \\
 T_{11}^n(r)=&\left(n^2-n-\frac{1}{2}k_\mathsf{T}^2r^2\right)j_n\left(k_\mathsf{L}r\right)+2k_\mathsf{L}r j_{n+1}\left(k_\mathsf{L}r\right), \\
 T_{13}^n(r)=&n\left(n+1\right)\left(\left(n-1\right)j_n\left(k_\mathsf{T}r\right)-k_\mathsf{T}rj_{n+1}\left(k_\mathsf{T}r\right)\right),\\
 T_{41}^{n}(r)=&\left(n-1\right)j_{n}\left(k_{\mathsf{L}}r\right)-k_{\mathsf{L}}r j_{n+1}\left(k_{\mathsf{L}}r\right), \\
 T_{42}^{n}(r)=&\frac{1}{2}r\left(\left(n-1\right)j_{n}\left(k_{\mathsf{T}}r\right)-k_{\mathsf{T}}r j_{n+1}\left(k_{\mathsf{T}}r\right) \right), \\
 T_{43}^n(r)=&\left(n^2-1-\frac{1}{2}k_\mathsf{T}^2 r^2\right)j_n\left(k_\mathsf{T}r\right)+k_\mathsf{T}r j_{n+1}\left(k_\mathsf{T}r\right), \\
 T_{51}^{n}(r)=&T_{41}^{n}(r),\\
 T_{52}^{n}(r)=&T_{42}^{n}(r),\\
 T_{53}^{n}(r)=&T_{43}^{n}(r),
\end{align}
where $j_n$ is the $n$-th spherical Bessel function, and $k_\mathsf{L}$ and $k_\mathsf{T}$ are the wave numbers for the longitudinal and the transverse wave, respectively, which have the following expressions:
\begin{align}
 k_\mathsf{L}=&\omega\sqrt{\frac{\rho_\mathsf{s}}{\lambda+2\mu}}, \\
 k_\mathsf{T}=&\omega\sqrt{\frac{\rho_\mathsf{s}}{\mu}}.
\end{align}
Also, $a_n^m$, $b_n^m$, $c_n^m$ and $d_n^m\in\mathbb{C}$ are coefficients of the spherical expansion.
Substituting \eqref{hatsigmartheta}, \eqref{hatsigmarphi} into \eqref{kyuten3} gives the following relation for $a_n^m$, $b_n^m$ and $c_n^m$:
\begin{align}
 b_n^m=&0, \label{bnm}\\
 T_{41}^n(\varepsilon) a_n^m+\frac{T_{43}^n(\varepsilon)}{k_\mathsf{T}}c_n^m=&0, \label{T41anm}
\end{align}
by which Eqs.~\eqref{kyutenhatur} and \eqref{kyutenhatsigmaij} are reduced as
\begin{align}
\hat{u}_{r}\left(\boldsymbol{x}\right)
&=\sum_{n=0}^{\infty}\sum_{m=-n}^{n}
\frac{a_{n}^{m}}{r}\left(U_{1}^{n}(r)-\frac{T_{41}^{n}(\varepsilon)}{T_{43}^{n}(\varepsilon)}U_{3}^{n}(r)\right)Y_{n}^{m}\left(\theta,\,\phi \right),
\label{kyutenhatur2} \\
 \hat{\sigma}_{rr}\left(\boldsymbol{x}\right)&=\sum_{n=0}^{\infty}\sum_{m=-n}^{n}\frac{2\mu}{r^2}a_n^m\left(T_{11}^n(r)-\frac{T_{41}^n(\varepsilon)}{T_{43}^n(\varepsilon)}T_{13}^n(r)\right)Y_n^m\left(\theta,\,\phi\right),
\label{kyutenhatsigrr2}
\end{align}
respectively.
By substituting \eqref{kyutenq}, \eqref{kyutendelq} and \eqref{kyutenhatur2} into \eqref{kyuten1}, 
substituting \eqref{kyutenp}, \eqref{kyutendelp} and \eqref{kyutenhatsigrr2} into \eqref{kyuten2}, 
and exploiting the orthogonal property of the spherical harmonics, 
we obtain a system of algebraic equations to determine $d_n^m$ and $a_n^m$, 
from which we obtain the asymptotic expansions of $\hat{u}_{r}$ and $\hat{\sigma}_{rr}$ as 
\begin{align}
 \hat{u}_{r}\left(\boldsymbol{x}\right)=&\frac{3}{2\mu k^2_{\mathsf{T}}+\rho_{\mathsf{f}}\omega^2}p_{,j}^0n_j\left(\boldsymbol{x}\right)
+\frac{k^2_\mathsf{L}}{\mu\left(4k_\mathsf{L}^2-3k_\mathsf{T}^2\right)}p^0\varepsilon+o\left(\varepsilon\right), \label{uhat_final} \\
 \hat{\sigma}_{rr}\left(\boldsymbol{x}\right)=&-p^0
+\left(
\frac{p^0_{,jj}}{9}
-\frac{3}{2\mu k_{\mathsf{T}}^2+\rho_\mathsf{f}\omega^2}p_{,j}^0n_j(\boldsymbol{x})
-\frac{p_{,ij}^0}{3}n_i(\boldsymbol{x})n_j(\boldsymbol{x})
\right)
\varepsilon+o\left(\varepsilon\right), \label{sigrrhat_final}
\end{align}
respectively. The asymptotic expansion of the traction $\hat{t}_i$ can easily be calculated from Eqs.~\eqref{kyuten3} and \eqref{sigrrhat_final}.
With these observations, we can evaluate the asymptotic expansion of $\delta J$ as
\begin{align}
 \delta J=&{\Re}\left[
\frac{4}{3}\pi\varepsilon^3\left(
\frac{3\left(\rho_\mathsf{s}-\rho_\mathsf{f}\right)}{2\rho_\mathsf{s}+\rho_\mathsf{f}}
p_{,j}\left(\boldsymbol{x}^0\right)\tilde{p}_{,j}\left(\boldsymbol{x}^0\right)
-
\frac{\Lambda_\mathsf{s}-\Lambda_\mathsf{f}}{\Lambda_\mathsf{s}\Lambda_\mathsf{f}}\rho_\mathsf{f}\omega^2
p\left(\boldsymbol{x}^0\right)\tilde{p}\left(\boldsymbol{x}^0\right)
\right)
\right] +o\left(\varepsilon^3\right),
\label{delJ}
\end{align}
where $\Lambda_\mathsf{s}$ is the bulk modulus for elastic material $\Omega^\mathsf{c}$ defined as
\begin{align}
 \Lambda_\mathsf{s}=\lambda+\frac{2}{3}\mu.
\end{align}

The topological derivative is defined as the coefficient of the leading term of the asymptotic expansion \eqref{delJ} of the objective function $J$ \cite{isakari2014topology,bonnet2007fm,novotny2003topological,carpio2008solving} as follows:
\begin{align}
 \delta
 J=\mathcal{T}_\Omega\left(\boldsymbol{x}\right)v\left(\varepsilon\right)
 +o\left(v\left(\varepsilon\right)\right),
 \label{deltaJ_tpodansei}
\end{align}
where $v\left(x\right)$ is a monotonically increasing function in $x>0$.
By comparing Eqs.~\eqref{delJ} and \eqref{deltaJ_tpodansei}, 
we obtain the topological derivative ${\mathcal T}_\Omega$ with respect to an appearance of an elastic material $\Omega_\varepsilon$ in the fluid matrix $\Omega$ as follows:
\begin{align}
 \mathcal{T}_\Omega \left(\boldsymbol{x}\right)
 ={\Re}\left[
\frac{3\left(\rho_\mathsf{s}-\rho_\mathsf{f}\right)}{2\rho_\mathsf{s}+\rho_\mathsf{f}}
p_{,j}\left(\boldsymbol{x}\right)\tilde{p}_{,j}\left(\boldsymbol{x}\right)
-
\frac{\Lambda_\mathsf{s}-\Lambda_\mathsf{f}}{\Lambda_\mathsf{s}\Lambda_\mathsf{f}}\rho_\mathsf{f}\omega^2
p\left(\boldsymbol{x}\right)\tilde{p}\left(\boldsymbol{x}\right)
\right],\label{danseitpoderi}
\end{align}
where $v(\varepsilon)$ in Eq.~\eqref{deltaJ_tpodansei} is chosen as $\displaystyle v\left(\varepsilon\right)=\frac{4}{3}\pi \varepsilon^3$.

The topological derivative ${\mathcal T}_{\Omega^\mathsf{c}}$ related to appearance of an acoustic material $\Omega_\varepsilon$ in the elastic inclusion $\Omega^\mathsf{c}$ can similarly be obtained follows:
\begin{align}
 \mathcal{T}_{\Omega^\mathsf{c}}\left(\boldsymbol{x}\right)=&
 {\Re}\left[\rho_\mathsf{f}\omega^2 \left(
 \left( A - B \right) \tilde{\sigma}_{ii}(\boldsymbol{x}) \sigma_{ii}(\boldsymbol{x})
 +3B\tilde{\sigma}_{ij}(\boldsymbol{x}) \sigma_{ij}(\boldsymbol{x})
 -\left(\rho_\mathsf{s}-\rho_\mathsf{f}\right)\omega^2\tilde{u}_i(\boldsymbol{x})u_i(\boldsymbol{x})
\right)
 \right], \label{72}
\end{align}
where the coefficients $A$ and $B$ are defined as
\begin{align}
 A=&\cfrac{3\left(\Lambda_{\mathsf{f}}-\Lambda_\mathsf{s}\right)\left(\lambda+2\mu\right)}
 {\left(3\lambda+2\mu\right)\left(12\mu\left(\lambda+\Lambda_\mathsf{f}\right)+9\lambda^2+4\mu^2\right)}, \\
 B=&\frac{5(\lambda+2\mu)}{2\mu\left(9\lambda+14\mu\right)},
\end{align}
respectively. 
Note that the result in \eqref{72} is consistent with the one in Guzina and Chikichev~\cite{guzina2007imaging}.

\subsection{A fast BEM-FEM solver for acoustic-elastic coupled problems}\label{bem_fem}
In order to evaluate the topological derivatives in Eqs.~\eqref{danseitpoderi} and \eqref{72}, 
we need to calculate the sound pressure $p$ and its gradient ${p_{,}}_j$ in $\Omega$, 
and the displacement $u_i$ and the stress $\sigma_{ij}$ in $\Omega^\mathsf{c}$, and their adjoint counterparts.
Although these quantities can appropriately be calculated by the boundary element method (BEM), a BEM-FEM (finite element method) coupled solver is utilised in this paper.
The proposed solver deals with the acoustic field and the elastic field by BEM and FEM, respectively.
This is because, in our future publications, 
we plan to extend the present topology optimisation for elastic materials other than the isotropic one, e.g. anisotropic material and Biot's poroelastic material for which the FEM is more suitable than the BEM.
In this section, we also present a fast algorithm for the BEM-FEM coupled solver, 
in which the fast multipole method (FMM) and a multi-frontal solver for sparse matrices are efficiently combined.

We present the formulation of the BEM-FEM coupled solver for the forward problem in Eqs.~\eqref{hel}--\eqref{rad}.
The adjoint problem in Eqs.~\eqref{danseiajohel}--\eqref{danseiajosommarfeld} can be solved in a same manner.

From Eqs.~\eqref{hel} and \eqref{bc_velo}, we have the following boundary integral equation:
\begin{align}
\frac{p\left(\boldsymbol{x}\right)}{2}&=\sum_{m=1}^{M^\mathsf{obs}}G(\boldsymbol{x}-\boldsymbol{x}_m^\mathsf{obs})
-\rho_\mathsf{f}\omega^2  \int_{\Gamma}
G\left(\boldsymbol{x}-\boldsymbol{y}\right)u_\ell\left(\boldsymbol{y}\right)n_\ell\left(\boldsymbol{y}\right)\mathrm{d}\Gamma\left(\boldsymbol{y}\right)
+\int_{\Gamma} \frac{\partial G\left(\boldsymbol{x}-\boldsymbol{y}\right)}{\partial n(\boldsymbol{y})}p\left(\boldsymbol{y}\right)\mathrm{d}\Gamma\left(\boldsymbol{y}\right),
 \label{danseibie2}
\end{align}
where $G(\boldsymbol{x})=e^{\mathsf{i}k_\mathsf{f}|\boldsymbol{x}|}/4\pi|\boldsymbol{x}|$ is the fundamental solution of three dimensional Helmholtz' equation.
The weak form in $\Omega^\mathsf{c}$ and the boundary condition in Eq.~\eqref{bc_power} gives the following equation:
\begin{align}
\int _\Gamma u^*_i n_ip \mathrm{d}\Gamma  
-\int _{\Omega^\mathsf{c}} u^*_{i,j}\sigma_{ji}\mathrm{d}\Omega  
+\rho_\mathrm{s}\omega^2\int _{\Omega^\mathsf{c}} u^*_{i}u_{i}\mathrm{d}\Omega 
=0,
\label{fem_2}
\end{align}
where $u_i^*$ is a test function. 
The proposed method solves the system of integral equations \eqref{danseibie2} and \eqref{fem_2}.
In the discretisation, $\Omega^\mathsf{c}$ is divided as $\Omega^\mathsf{c}=\cup_{e=1}^{N_{\mathsf{fe}}} \Omega_e$, 
where $\Omega_e$ is a tetrahedron, 
and $\Gamma$ is divided as $\Gamma=\cup_{j=1}^{N_\mathsf{be}} \Gamma_j$, 
where $\Gamma_j$ is a triangular patch which coincides to a surface of a tetrahedron $\Omega_e$.
In this study, the sound pressure $p$ is approximated by locally constant functions on $\Gamma_j$ 
while the displacement $u_i$ is approximated by locally linear functions in $\Omega_e$.
With these settings, a standard collocation for Eq.~\eqref{danseibie2} and the Galerkin discretisation, 
in which the test functions are also expanded as the locally linear functions in $\Omega_e$, 
for Eq.~\eqref{fem_2} gives the following system of algebraic equations:
\begin{align}
\begin{pmatrix}
\frac{\mathsf{I}}{2}-\mathsf{D} & \mathsf{S} \\
\mathsf{N} & \mathsf{K}-\rho_\mathsf{s}\omega^2\mathsf{M}
\end{pmatrix}
\begin{pmatrix}
\mathsf{p} \\
\mathsf{u}
\end{pmatrix}
=
\begin{pmatrix}
\mathsf{p}^{\mathsf{src}}\\
\mathsf{0}
\end{pmatrix},
\label{renritupu}
\end{align}
where $\mathsf{S}$ and $\mathsf{D}$ are the coefficient matrices for single and double layer potentials, respectively.
$\mathsf{K}$ and $\mathsf{M}$ are the finite element stiffness and mass matrices, respectively.
Also, the matrix $\mathsf{N}$ is stemmed from the first term of Eq.~\eqref{fem_2}.
The vectors $\mathsf{p}$ and $\mathsf{p}^\mathsf{src}$ contain the total and the incident sound pressures, respectively, on the collocation points $\boldsymbol{x}_i\in\Gamma_j$, 
and the vector $\mathsf{u}$ is composed of nodal displacements in $\Omega_e$.

In the following, we propose a fast solver for the algebraic equation \eqref{renritupu}.
The basic idea of the proposed solver is to combine the fast multipole method (FMM)~\cite{rokhlin1985rapid,greengard1987fast} for calculations of the BEM matrices,
and a multi-frontal solver to factorise the FEM matrix in by exploiting the block structure in Eq.~\eqref{renritupu}.
The second row in \eqref{renritupu} can be written as 
\begin{align}
\mathsf{u}=&(\mathsf{K}-\rho_\mathsf{s}\omega^2\mathsf{M})^{-1}\mathsf{N}\mathsf{p},
\label{kousokuu}
\end{align}
provided that $\omega^2$ is not an eigenvalue of a homogeneous Neumann problem in $\Omega^\mathsf{c}$.
Note that the matrix $\mathsf{K}-\rho_{\mathsf{s}}\omega^2\mathsf{M}$ is a sparse matrix which can efficiently be factorised by a multi-frontal solver.
By substituting Eq.~\eqref{kousokuu} into the first row of Eq.~\eqref{renritupu}, 
we obtain the following algebraic equation:
\begin{align}
\left(\frac{\mathsf{I}}{2}-\mathsf{D}+\mathsf{S}\left(\mathsf{K}-\rho_\mathsf{s}\omega^2\mathsf{M}\right)^{-1}\mathsf{N}\right)\mathsf{p}
=
\mathsf{p}^{\mathsf{src}}.
\label{GMRES}
\end{align}
From the Calderon identity, one observes that almost all the eigenvalues of the matrices $\mathsf{S}$ and $\mathsf{D}$ are close to 0~\cite{christiansen2002preconditioner,niino2012preconditioning,isakari2012calderon}.
Thus, the condition number of the coefficient matrix in Eq.~\eqref{GMRES} is expected to be small, 
which leads fast convergence of an iterative solver such as GMRES.
Matrix-vector products involved in an iterative solver can efficiently be performed by either the FMM or a multi-frontal solver.
The algorithm for solving the boundary value problem in Eqs.~\eqref{hel}--\eqref{rad} is summarised as follows:
\begin{enumerate}
 \item The matrix $\mathsf{K}-\rho_\mathsf{s}\omega^2\mathsf{M}$ in \eqref{kousokuu} is factorised by a multi-frontal solver, and the factorised matrices are stored.
 \item Equation \eqref{GMRES} is solved with an iterative solver to obtain the sound pressure $p$ on the collocation point $\boldsymbol{x}_j\in\Gamma_{j}~(j=1,...,N_\mathsf{be})$.
       The product of the coefficient matrix in Eq.~\eqref{GMRES} and a vector $\mathsf{x}$, which is required in the algorithm of the iterative solver, is performed as follows:
       \begin{enumerate}
	\item The matrix-vector product $\mathsf{z}:=\mathsf{N}\mathsf{x}$ is calculated. In the calculation, the sparsity of the matrix $\mathsf{N}$ is exploited.
	\item The vector $\mathsf{y}:=(\mathsf{K}-\rho_\mathsf{s}\omega^2\mathsf{M})^{-1}\mathsf{z}$ is calculated by solving the algebraic equation $(\mathsf{K}-\rho_\mathsf{s}\omega^2\mathsf{M})\mathsf{y}=\mathsf{z}$ with a multi-frontal solver.
	      Note that the coefficient matrix $\mathsf{K}-\rho_\mathsf{s}\omega^2\mathsf{M}$ has already been factorised in 1.
	\item $\left(\frac{\mathsf{I}}{2}-\mathsf{D}\right)\mathsf{x+Sy}$ is calculated by the FMM.
       \end{enumerate}
 \item The matrix-vector product $\mathsf{z}':=\mathsf{Np}$ is calculated. In the calculation, the sparsity of the matrix $\mathsf{N}$ is, again, exploited.
 \item The displacement $u$ on the nodal point of the finite elements $\Omega_e$ is obtained by solving $(\mathsf{K}-\rho_{\mathsf{s}}\omega^2\mathsf{M})\mathsf{u}=\mathsf{z}'$.
       Again, note that the coefficient matrix $\mathsf{K}-\rho_{\mathsf{s}}\omega^2\mathsf{M}$ has already been factorised in 1.
 \item The sound pressure $p$ and its gradient $p_{,i}$ are calculated by integral representations on arbitrary points $x\in\Omega$.
 \item The displacement $u$ and the stress $\sigma_{ij}$ are calculated by interpolation with shape functions on arbitrary points $x\in\Omega^\mathsf{c}$.
\end{enumerate}

In the implementation of the FMM, a low-frequency FMM \cite{isakari2014topology} is employed.

\subsection{A level-set-based topology optimisation}\label{lsm}
In this subsection, we briefly review a level-set-based methodology to solve the optimisation problem 
to find an optimal distribution of elastic material(s) $\Omega^\mathsf{c}\subset D$ which minimises the objective function in Eq.~\eqref{obj} subject to the constrain conditions in Eqs.~\eqref{hel}--\eqref{rad} (see also Figure~\ref{01formulation_optimisation}).
The reader is referred to the original paper~\cite{yamada2010topology} and our previous papers~\cite{isakari2016multi,isakari2014topology,jing2015level} for further details.

In the level set method, domains $\Omega$ and $\Omega^\mathsf{c}$, and its boundary $\Gamma$ is recognised as 
\begin{align}
 \Omega^\mathsf{c}&=\{\boldsymbol{\xi}~|~0<\phi(\boldsymbol{\xi})\le 1\},\label{phi1} \\  
 \Gamma           &=\{\boldsymbol{\xi}~|~\phi(\boldsymbol{\xi})=0\}, \label{phi2}\\ 
 \Omega           &=\{\boldsymbol{\xi}~|~-1<\phi(\boldsymbol{\xi})\le 0\}, \label{phi3}
\end{align}
respectively.
With the level set function $\phi$, the topology optimisation problem is converted into the problem to find an optimal distribution of $\phi$ in the design domain $D$
which minimises the objective function in Eq.~\eqref{obj} under the constrains in Eqs.~\eqref{hel}--\eqref{rad}.
We explore the optimum distribution of $\phi$ using the topological derivative ${\mathcal T}_\Omega$ and ${\mathcal T}_{\Omega^\mathsf{c}}$ in Eqs.~\eqref{danseitpoderi} and \eqref{72} from an initial distribution $\phi_0(\boldsymbol{\xi})$ as follows:
\begin{align}
 \frac{\partial\phi(\boldsymbol{\xi},t)}{\partial t}=-C{\mathcal T}_\Omega(\boldsymbol{\xi},t)+\tau L^2\nabla^2 \phi(\boldsymbol{\xi}) \ \ \ \ \ \mathrm{for} \ \ \phi(\boldsymbol{\xi},t)<0, \label{timeeq1}\\
 \frac{\partial\phi(\boldsymbol{\xi},t)}{\partial t}=C{\mathcal T}_{\Omega^\mathsf{c}}(\boldsymbol{\xi},t)+\tau L^2\nabla^2 \phi(\boldsymbol{\xi}) \ \ \ \ \ \mathrm{for} \ \ \phi(\boldsymbol{\xi},t)>0, \label{timeeq2} \\
 \phi(\boldsymbol{\xi},0)=\phi_0(\boldsymbol{\xi}), \label{time_initial}
\end{align}
where $t$ represents a fictitious time, $C>0$ is a constant, and $L$ is a characteristic length of the design domain $D$.
In the present method, the topological derivatives is used to modify the distribution of the level set function $\phi$.
${\mathcal T}_\Omega(\boldsymbol{\xi}, t)$ in Eq.~\eqref{timeeq1}, for example, works to allocate a small elastic scatterer on $\boldsymbol{\xi}$ when ${\mathcal T}_\Omega$ is negative by increasing $\phi(\boldsymbol{\xi}, t)$.
Also, $\tau>0$ is a parameter which prescribes the complexity of the geometry of $\Omega^\mathsf{c}$~\cite{yamada2010topology}.
The following boundary condition for the time evolution equations \eqref{timeeq1} and \eqref{timeeq2} is also defined: 
\begin{align}
\phi(\boldsymbol{\xi},t) =\bar{\phi(\boldsymbol{\xi},t)}<0 \ \ \ \ \mathrm{for} \ \ \boldsymbol{\xi}\in \partial D \ \mathrm{and} \ t>0,
\label{timebou}
\end{align}
where $\bar{\phi}$ is a known function.
The boundary condition \eqref{timebou} is imposed so that $\Omega^\mathsf{c}\subset D$ holds.

Thus, the optimisation problem is now converted to the initial-boundary value problem \eqref{timeeq1}--\eqref{timebou}, 
which is solved with FEM in this study.

\subsection{Algorithm of the present topology optimisation} \label{tpo_algo}
Combining all the techniques presented above, the algorithm of the proposed topology optimisation is summarised as follows:
\begin{enumerate}
 \item The fixed design domain $D$ is divided into finite elements (voxels) as $D=\cup_{e=1}^{N_\mathsf{D}}\Omega^\mathsf{D}_e$, where $\Omega^\mathsf{D}_e$ is a voxel, and $N_\mathsf{D}$ is the number of the voxels.
 \item An initial distribution of the level set function $\phi(\boldsymbol{\xi},0)$ is given on nodes (lattice points) of the finite elements $\cup_{e=1}^{N_\mathsf{D}}\Omega^\mathsf{D}_e$.
 \item A set $X$ of the points $\boldsymbol{x}$ such that $\phi(\boldsymbol{x})=0$ on the lattice edge is explored.
 \item By appropriately connecting the elements of $X$, triangular patches that cover the iso-surface of zero value of the level set function, and the triangular patches are stored in STL format.
       For details on this procedure, the reader is referred to \cite{shichi2012level}.
 \item The STL data is converted to the boundary elements $\Gamma=\cup_{j=1}^{N_\mathsf{be}} \Gamma_j$ and the finite elements $\Omega^\mathsf{c}=\cup_{e=1}^{N_{\mathsf{fe}}} \Omega_e$. 
       For the meshing, NETGEN~\cite{netgen} is used in this study.
       Note that NETGEN can not only generate boundary/finite meshes but also improve the quality of the meshes.
 \item The boundary value problem \eqref{hel}--\eqref{rad} (the forward problem) defined in $\overline{\Omega\cup\Omega^\mathsf{c}}$ is solved by the BEM-FEM coupled solver presented in Section \ref{bem_fem}, 
       and evaluate the objective function in Eq.~\eqref{obj}. When the objective function converges, stop.
 \item The boundary value problem \eqref{danseiajohel}--\eqref{danseiajosommarfeld} (the adjoint problem) defined in $\overline{\Omega\cup\Omega^\mathsf{c}}$ is solved by the BEM-FEM coupled solver in Section \ref{bem_fem}.
       Note that the procedure 1 in the algorithm of the BEM-FEM coupled solver can be skipped since matrix $\mathsf{K}-\rho_\mathsf{s}\omega^2\mathsf{M}$ has already been factorised in the forward analysis.
       Note also that some of the quantities in the FMM algorithms such as tree structures, direct interactions and M2L operators, etc calculated in the forward analysis can be recycled in the adjoint procedure.
 \item With the state  and adjoint variables calculated in 6 and 7, the topological derivatives in Eqs.~\eqref{danseitpoderi} and \eqref{72} on all the lattice points expanding $D$ are evaluated.
 \item The initial-boundary value problem \eqref{timeeq1}--\eqref{timebou} is solved by FEM. Go to 3.
\end{enumerate}

\section{Numerical examples}
In this section, we present some numerical examples to confirm the validity and the efficiency of the proposed method. 
We first state common issues to all the examples to follow:
\begin{itemize}
 \item \verb@PARDISO@ routines by Intel MKL is used as the multi-frontal solver involved in the present BEM-FEM solver.
 \item GMRES is used as the iterative solver involved in the present BEM-FEM solver.
 \item The tolerance of GMRES in the present BEM-FEM solver is set to be $10^{-5}$.
 \item Truncation numbers for infinite series in the FMM are numerically determined so that the truncation error is less than $10^{-5}$.
 \item All numerical experiments were run on a PC with Intel Xeon CPU E5-4650 with 32 cores. The code is OpenMP parallelised.
\end{itemize}

\subsection{Performance tests of the present BEM-FEM coupled solver}\label{performance_test}
In this subsection, we check the performance of the proposed BEM-FEM coupled solver presented in Section \ref{bem_fem} with benchmark problems.
To this end, the algebraic equations \eqref{renritupu} corresponding to benchmark problems are solved by conventional direct and iterative solvers as well as by the present solver.
We used a LAPACK routine \verb@ZGESV@ and GMRES as the conventional direct and iterative solvers to naively solve Eq.~\eqref{renritupu}, respectively.
In the conventional GMRES, the FMM is not employed to accelerate the matrix-vector products involved in the algorithm of GMRES, and the tolerance of the GMRES is set to be $10^{-12}$, 
which is set by numerical experiments so that the total errors for the sound pressures and the displacements are comparable to those by the present method and the conventional direct solver.

As the first benchmark problem, we consider a sound scattering by an elastic sphere whose analytical solution is available in \cite{eringen1978elastodynamics}.
We set an elastic sphere $\Omega^\mathsf{c}:=\{\boldsymbol{x}~|~|\boldsymbol{x}|< 0.25\}$ in an acoustic host matrix.
We here assume that the elastic inclusion and host acoustic matrix are composed of a tungsten and water, respectively.
The parameters for the tungsten are set as the density $\rho_\mathsf{s}=64.85$, Young's modulus $E=174.57$ and Poisson's ratio $\nu=0.3$, 
which are normalised by the material parameters of water as the density $\rho_\mathsf{f}=1.0$ and the bulk modulus $\Lambda_\mathsf{f}=1.0$.
As the incident wave, we used a plane wave propagating in $x_3$ direction with the frequency $\omega=1.0$.
The amplitude of the incident wave is set to be 1.0.

\begin{figure}[h]
 \begin{center}
  \includegraphics[scale=0.5]{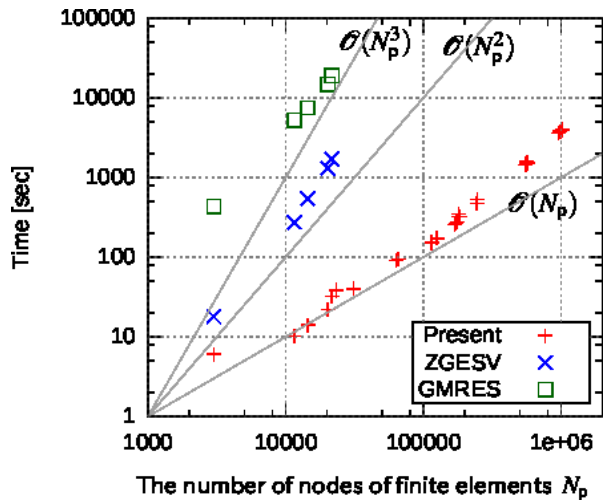}
  \caption{Computational time for the present and conventional solvers against the number of finite element nodes $N_\mathsf{p}$ for a sound scattering problem by an elastic sphere.}
  \label{04time}
 \end{center}
\end{figure}

With these settings, the algebraic equations \eqref{renritupu} are solved by the present method, 
in which the surface of $\Omega^\mathsf{c}$ is divided into $2456$ boundary elements, 
and $\Omega^\mathsf{c}$ into finite elements of $3007$ nodal points.
The average of the relative errors for the sound pressures on collocation points and 
for nodal values of the displacements were 0.0628\% and 0.107\%, respectively.
The accuracy of the present method is comparable to that of a conventional direct solver and a conventional iterative solver.
We then discuss the timing. 
Figure~\ref{04time} shows the computational time for the present and conventional solvers against the number of the nodes of finite elements $N_\mathsf{p}$.
The computational time is measured using an OpenMP run-time library routine \verb@omp_get_wtime@.
One observes that the computational time for the present method scales as $N_\mathsf{p}$, 
and the present method is the fastest among the tested solvers even when the degrees of freedom is relatively small.
The computational complexities of the conventional iterative and direct solvers are ${\mathcal O}(N_\mathsf{p}^2)$ and ${\mathcal O}(N_\mathsf{p}^3)$, respectively.
We think that the bad condition of the coefficient matrix in Eq.~\eqref{renritupu} makes the convergence of the conventional iterative solver slow.
Indeed, the iteration number for the conventional GMRES is 390 in the case of $N_\mathsf{p}=3007$.
On the other hand, the condition of the coefficient matrix in Eq.~\eqref{GMRES} in the present method is well as discussed in Section \ref{bem_fem}.
The iteration number for Eq.~\eqref{GMRES} is 3 for this case.

Table~\ref{table2} shows computational time for each procedure in the present BEM-FEM solver in the case of $N_\mathsf{p}=194536$.
One finds that calculating the FEM matrices $\mathsf{N}$ and $\mathsf{K}-\rho_\mathsf{s}\omega^2\mathsf{M}$, and factorising $\mathsf{K}-\rho_\mathsf{s}\omega^2\mathsf{M}$ take almost half of the whole computational time.
As indicated in Section \ref{tpo_algo}, after solving the forward problem in \eqref{hel}--\eqref{rad} in the process of the optimisation, 
we need not to repeat these procedures in solving the adjoint problem in Eq.~\eqref{danseiajohel}--\eqref{danseiajosommarfeld}
since the FEM matrices $\mathsf{N}$ and $\mathsf{K}-\rho_\mathsf{s}\omega^2\mathsf{M}$ for the adjoint problem are common to those for the forward problem.
With these observations, it is expected that the computational time for solving the adjoint problem is as approximately half as that for solving the forward problem.
Thus, the present BEM-FEM solver can efficiently be applied to the topology optimisation.
\begin{table}[h]
\begin{center}
\caption{Computational time for each procedure in the present BEM-FEM solver for sound scattering problem by an elastic sphere in the case of $N_\mathsf{p}=194536$.}
\label{table2}
 \begin{tabular}{c|r}
  Procedure & Comp. time [sec]\\
  \hline
  Calculation of matrices $\mathsf{N}$ and $\mathsf{K}-\rho_\mathsf{s}\omega^2\mathsf{M}$ & 24.44 \\
  Factorisation of $\mathsf{K}-\rho_\mathsf{s}\omega^2\mathsf{M}$ by \verb@PARDISO@ & 51.78 \\
  GMRES for solving Eq.~\eqref{GMRES} & 70.36 \\
  Othres & 6.05 \\
  \hline
  Sum & 152.63
 \end{tabular}
 \end{center}
  \end{table}

\begin{figure}[h]
 \begin{center}
  \includegraphics[scale=0.35]{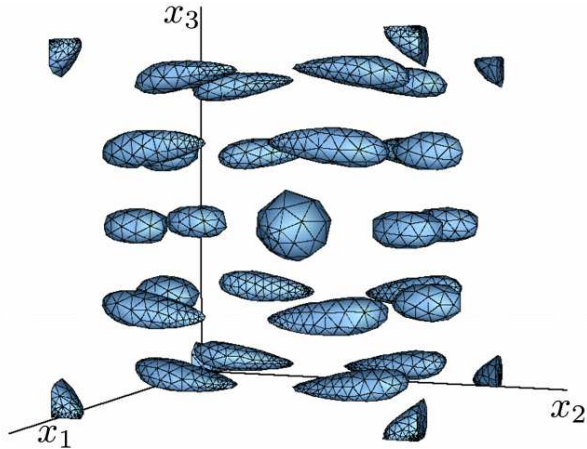}
  \caption{A complex-shaped elastic scatterers taken from \cite{isakari2014topology}.}
  \label{05fukuzatu}
 \end{center}
\end{figure}

As the second benchmark problem, we consider elastic scatterers shown in Figure~\ref{05fukuzatu} to check the applicability of the present BEM-FEM solver to a complex-shaped domain.
The elastic scatterers in Figure~\ref{05fukuzatu} are taken from \cite{isakari2014topology}, 
which are obtained in a process of a topology optimisation.
The elastic scatterers are expressed with finite elements of $N_\mathsf{p}=3845$.
The material parameters and the incident wave are set as in the previous test.
In this benchmark problem, we compare the performance of the present solver with that of a conventional iterative solver; GMRES is naively employed to solve the algebraic equation \eqref{renritupu}.
Table~\ref{table3} shows the computational time for scattering by the complex-shaped scatterers and by a sphere with $N_\mathsf{p}=3007$ for comparison.
The computational time of the present solver for the complicated shape is as only 9 times as that for the sphere, 
while the conventional iterative solver for the complicated shape is much slower than the case for a sphere.
Thus, the performance of the present solver is less affected by shape of scatterers than that of the conventional iterative solver.
This fact can also be seen in Table~\ref{table4}, which shows the iterative number of the GMRES.
We confirm that, for complicated domain, the condition number of Eq.~\eqref{renritupu} can be quite large, while that of Eq.~\eqref{GMRES} is kept relatively small.
This is because the present method exploits the spectral properties of the boundary element matrices $\mathsf{S}$ and $\mathsf{D}$ according to the Calderon formulae.
Thus, the proposed BEM-FEM solver is efficient for complex-shaped domain which is often encountered in topology optimisation.
\begin{table}[h]
 \begin{center}
  \caption{Computational time of the present and a conventional iterative solver for sound scattering problems by elastic scatterers.}
  \label{table3}
  \begin{tabular}{c|rr}
   Method & Sphere & Complicated shape \\
   \hline
   Present & 6.04 sec & 52.87 sec \\
   Conventional GMRES  & 421.44 sec & 2147659709.14 sec
  \end{tabular}
  \caption{The number of iteration for the present and a conventional iterative solver for sound scattering problems by elastic scatterers.}
  \label{table4}
  \begin{tabular}{c|rr}
   Method & Sphere & Complicated shape \\
   \hline
   Present & 3 & 74 \\
   Conventional GMRES & 390 & 16092
  \end{tabular}
 \end{center}
\end{table}

\subsection{Verification of the topological derivative}
\begin{figure}[h]
 \begin{center}
  \includegraphics[scale=0.6]{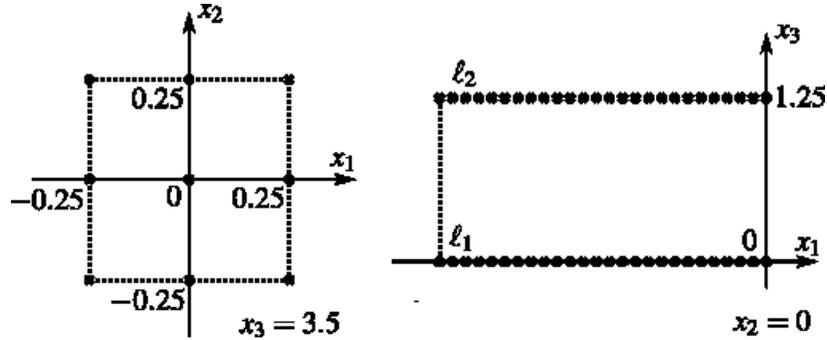}
  \caption{Definitions for (left:) the Observation points in \eqref{ex:obj} and (right:) the evaluation points for the topological derivative for the verification of the topological derivative in Eq.~\eqref{danseitpoderi}.}
  \label{15obj_and_td}
 \end{center}
\end{figure}

\begin{figure}[h]
 \begin{center}
  \includegraphics[scale=0.4]{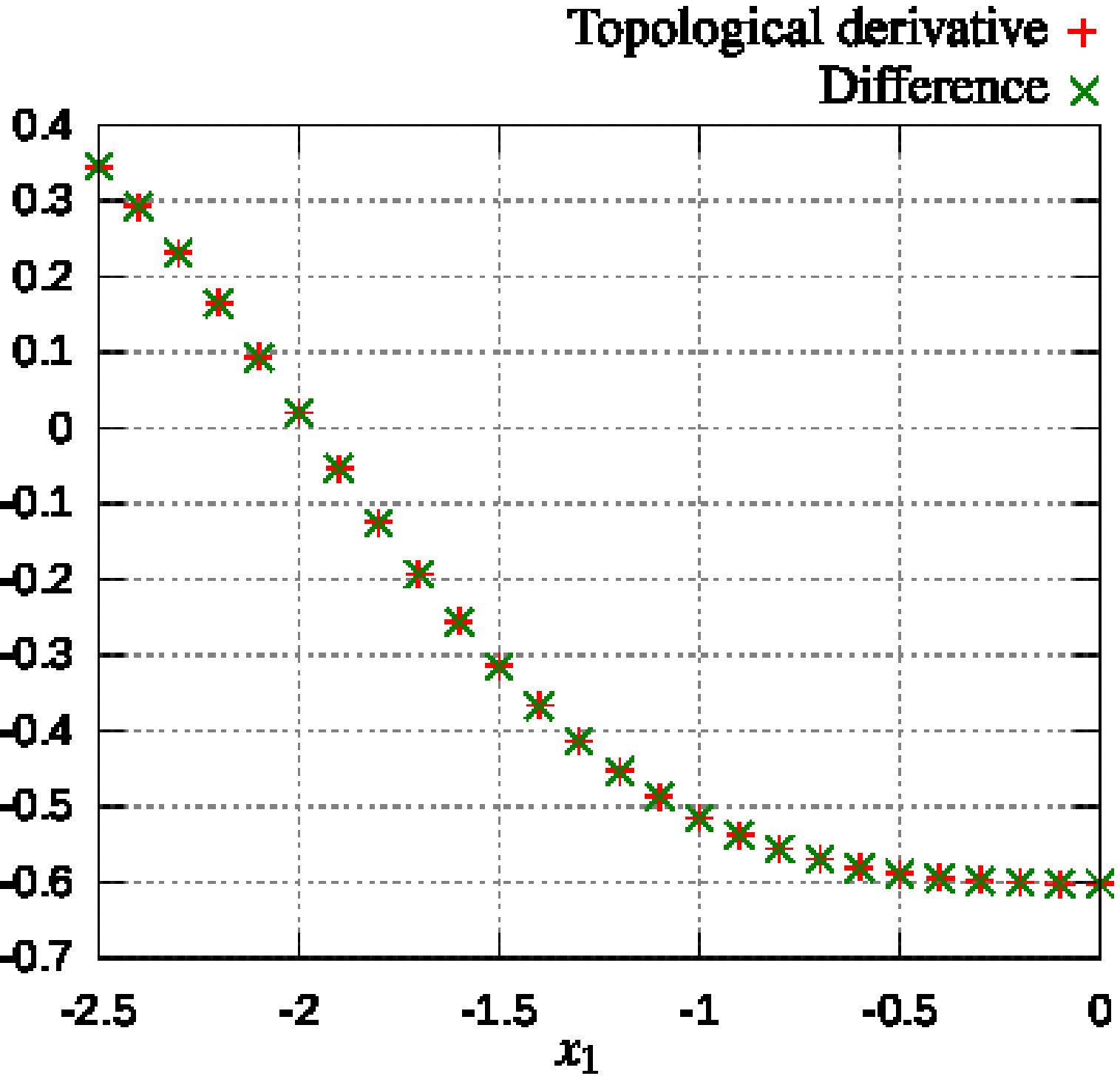}
  \caption{Topological derivatives on the line $\ell_1$ in Figure~\ref{15obj_and_td}.}\vspace{15pt}
  \label{06sabun}
  \includegraphics[scale=0.4]{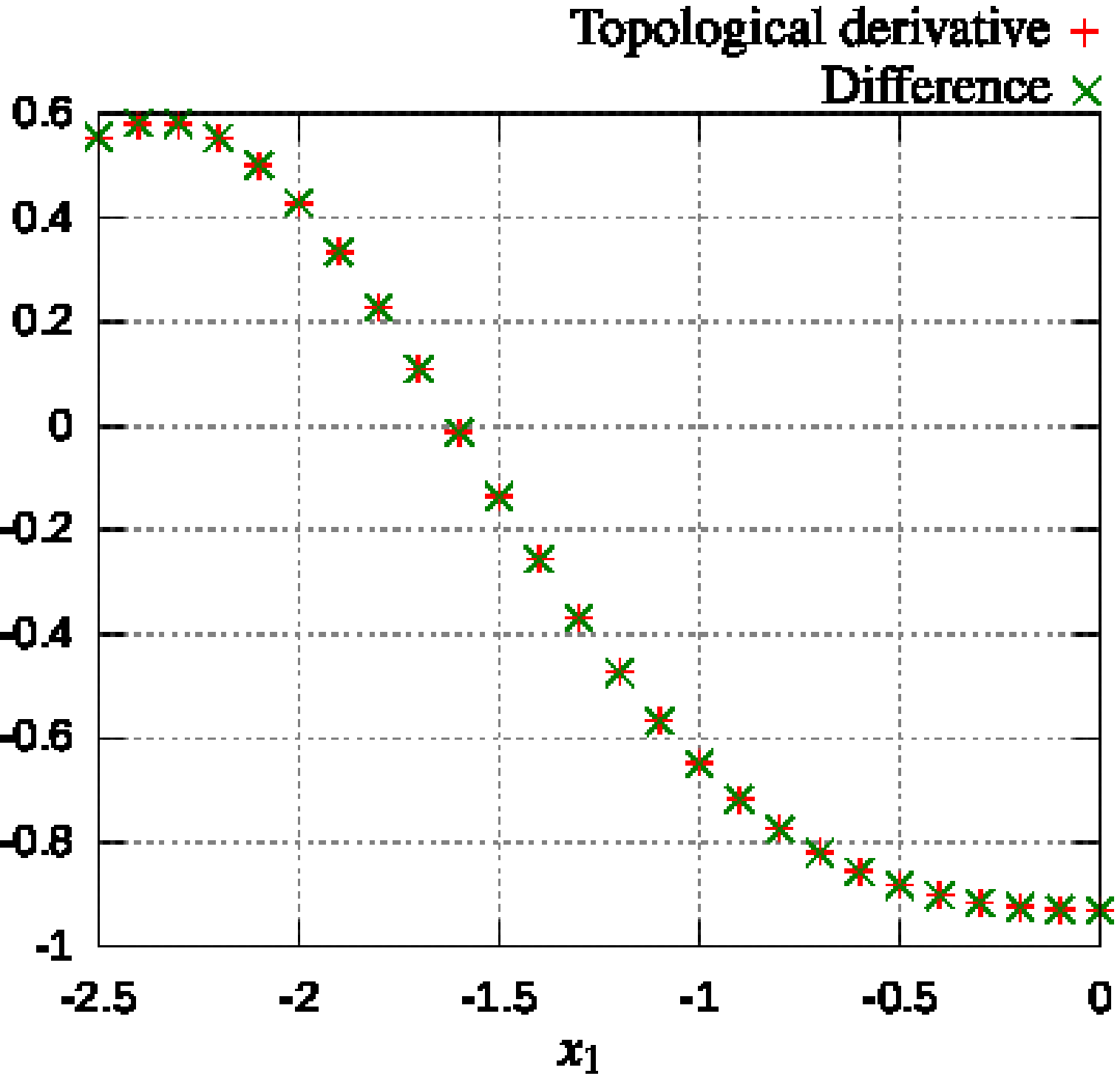}
  \caption{Topological derivatives on the line $\ell_2$ in Figure~\ref{15obj_and_td}.}
  \label{07sabun}
 \end{center}
\end{figure}

In this subsection, we numerically verify the topological derivative ${\mathcal T}_\Omega$ in Eq.~\eqref{danseitpoderi} derived in Section \ref{sec:td}.
We here consider to put an spherical elastic material of infinitesimal radius in $\mathbb{R}^3$ filled with an acoustic host matrix.
We assume that the elastic sphere and the acoustic host matrix are composed of acrylonitrile butadiene styrene (ABS) resin and water, respectively.
The material parameters for ABS resin is normalised as the density $\rho_\mathsf{s}=1.1$, the Young modulus $E=1.17$ and the Poisson ratio $\nu=0.369$ with the ones for water as the density $\rho_\mathsf{f}=1.0$ and the bulk modulus $\Lambda_\mathsf{f}=1.0$.
The incident wave is assumed to be a plane wave propagating in $x_3$ direction with the frequency $\omega=1.0$.
The amplitude of the incident wave is set as 1.0.
We define the following sum of the sound norm on the 9 observation points $\boldsymbol{x}_m^\mathsf{obs}$ (defined in the left figure of Figure~\ref{15obj_and_td}) as the objective functional:
\begin{align}
J=\frac{1}{2}\sum_{m=1}^{9}|p(\boldsymbol{x}_m^\mathsf{obs})|^2.
\label{ex:obj} 
\end{align}
Figures \ref{06sabun} and \ref{07sabun} show the distribution of the topological derivatives for the objective functional in Eq.~\eqref{ex:obj} on lines $\ell_1$ and $\ell_2$ in the right figure of Figure~\ref{15obj_and_td}, respectively.
In the figures, the ``topological differences'' $\mathcal D$ are also plotted, which are defined as follows:
\begin{align}
 {\mathcal D}=\frac{J_{\Omega\setminus\overline{\Omega_\varepsilon}}-J_\Omega}{v(\varepsilon)},
 \label{ex:sabun}
\end{align}
where $v(\varepsilon)=4\pi\varepsilon^3/3$, 
and $J_\Omega$ and $J_{\Omega\setminus\overline{\Omega_\varepsilon}}$ represent the objective function before and after a small spherical elastic scatterer of radius $\varepsilon$ is introduced, respectively.
In the case that $\varepsilon$ is small in Eq.~\eqref{ex:sabun}, $\mathcal D$ is expected to agree with the topological derivative.
We used the BEM for the calculation of $J_\Omega$ and $J_{\Omega\setminus\overline{\Omega_\varepsilon}}$ with $\varepsilon=0.03$ whose surface is divided into 2000 boundary elements.
One confirms that the topological derivatives derived in this paper agree well with the reference values.

\subsection{Optimal designs}
In this subsection, we show two examples of optimal design of sound scatterers which reduces the sound norm on some preset observation points.
In the first example, we consider a hard elastic material which may appropriately be modelled by a rigid one.
We show the obtained configuration of the hard elastic material is similar to that of rigid one, 
with which the validity of the proposed topology optimisation is confirmed.
In the second example, we consider a design problem of a soft elastic sound scatterer which cannot be solved by conventional topology optimisation methods.

\subsubsection{Hard scatter}
\begin{figure}[h]
  \begin{center}
   \includegraphics[scale=0.15]{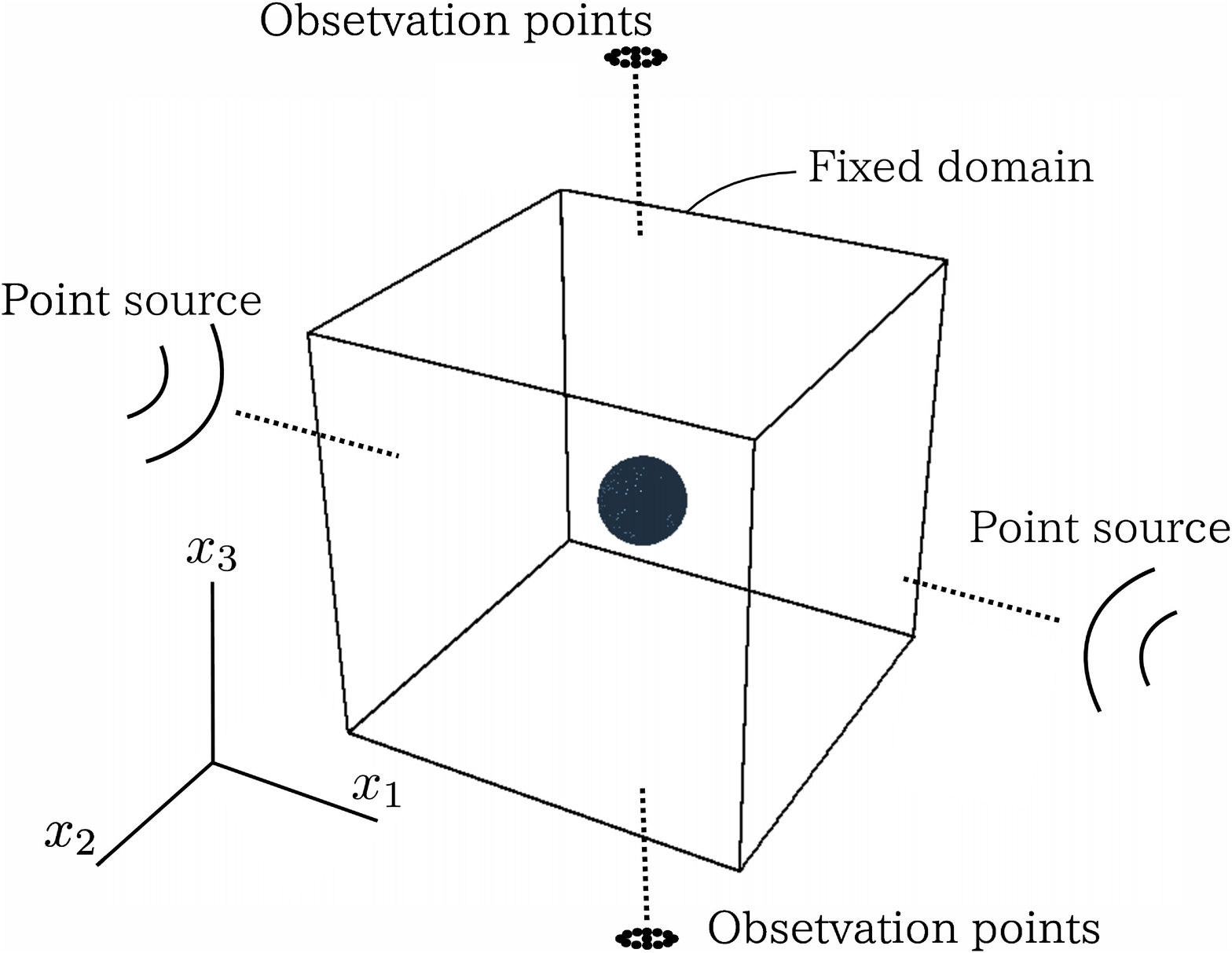}
   \caption{Settings for a topological optimal design for sound scatterer with hard elastic material. In the figure, the initial configuration of the elastic material is also plotted.}
   \label{075zyouken}
  \end{center}
  \end{figure}

We explore, with the present topology optimisation, an optimal distribution of hard elastic material in a design domain $D:=\{\boldsymbol{x}~|~0\le x_i \le 2.5 (i=1,2,3)\}$, 
which minimises the sound norm on observation points.
The observation points are set as $(1.25, 1.25, 5.0)$ and $(1.25, 1.25, -2.5)$, and 20 points on the circles, whose centre are these points and radius is $1.0$.
The circles are parallel to $x_1x_2$ plane (Figure~\ref{075zyouken}).
Sound sources are set on $(5.0, 1.25, 1.25)$ and $(-2.5, 1.25, 1.25)$ whose intensity and frequency are set as 70 and $2\pi$, respectively.
As the initial guess for the elastic material, we used a sphere whose centre and radius are $(1.25, 1.25, 1.25)$ and $0.25$, respectively.
The elastic material and the host acoustic material are respectively assumed to be a tungsten and water whose material parameters are listed in Section \ref{performance_test}.
In the optimisation algorithm in Section \ref{lsm}, we divide the design domain $D$ into $100\times 100\times 100$ finite elements, and set $\tau$ in Eqs.~\eqref{timeeq1} and \eqref{timeeq2} as $10^{-4}$.

\begin{figure}[h]
  \begin{center}
   \includegraphics[scale=0.2]{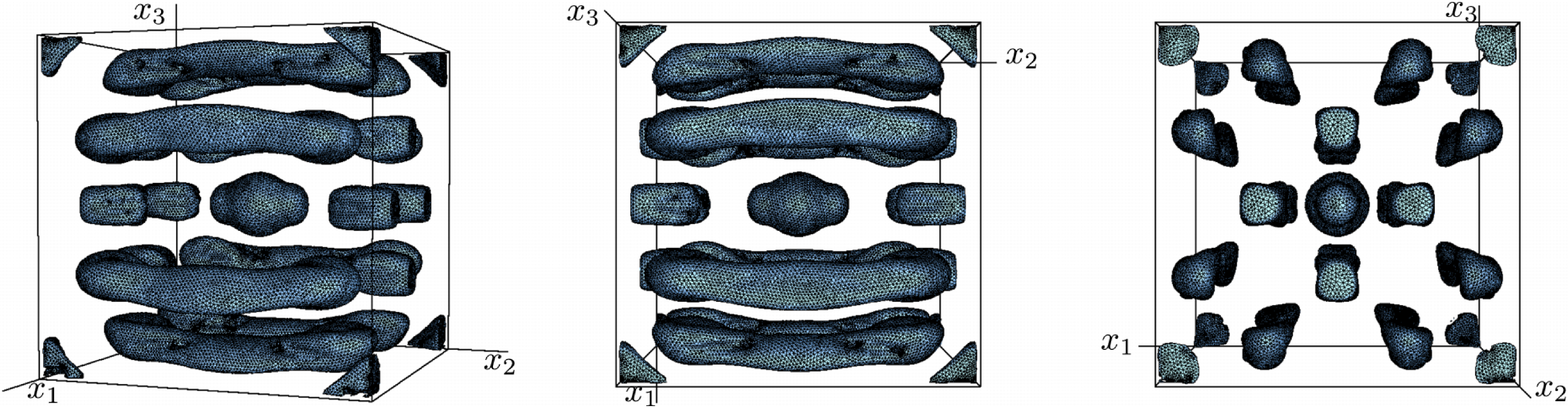}
   \caption{The obtained configuration of tungsten.}
   \label{08tang_mesh}
   \includegraphics[scale=0.2]{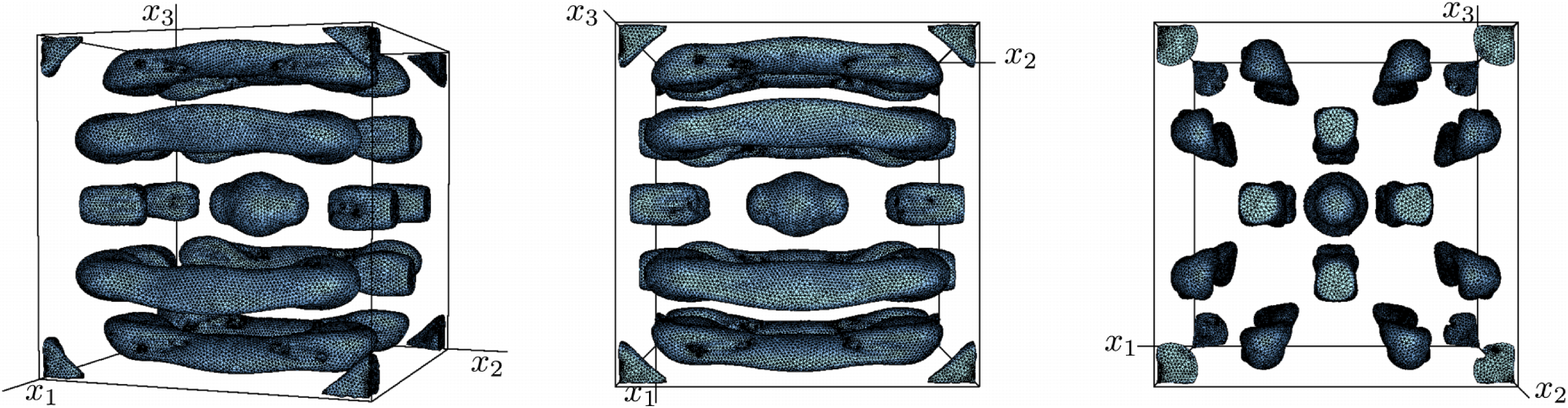}
   \caption{The obtained configuration of rigid material.}
   \label{09rigid_mesh}
  \end{center}
 \end{figure}

We show the obtained configuration in Figure~\ref{08tang_mesh}.
For comparison, we also show in Figure~\ref{09rigid_mesh} a result of the optimisation problem in which the elastic material is replaced by a rigid one.
One finds that the obtained configuration of tungsten is similar to that of rigid material, 
which is reasonable since the tungsten is quite hard and heavy compared to water.
Thus, the tungsten embedded in water can appropriately be approximated by rigid material.
This fact can also be confirmed by the expression of the topological derivative in \eqref{danseitpoderi}.
By taking limits as $\rho_\mathsf{f}/\rho_\mathsf{s}\rightarrow 0$ and $\Lambda_\mathsf{f}/\Lambda_\mathsf{s}\rightarrow 0$ (sound-hard limit), 
the topological derivative ${\mathcal T}_\Omega$ becomes as
\begin{align}
{\mathcal T}_\Omega(\boldsymbol{x}) \rightarrow \Re\left[\frac{3}{2}p_{,j}(\boldsymbol{x})\tilde{p}_{,j}(\boldsymbol{x})-\rho_\mathsf{f}\omega^2p(\boldsymbol{x})\tilde{p}(\boldsymbol{x}) \right],
\end{align}
which is identical to the topological derivative for the rigid material \cite{isakari2014topology,bonnet2007fm}.

\begin{figure}[h]
  \begin{center}
   \includegraphics[scale=0.35]{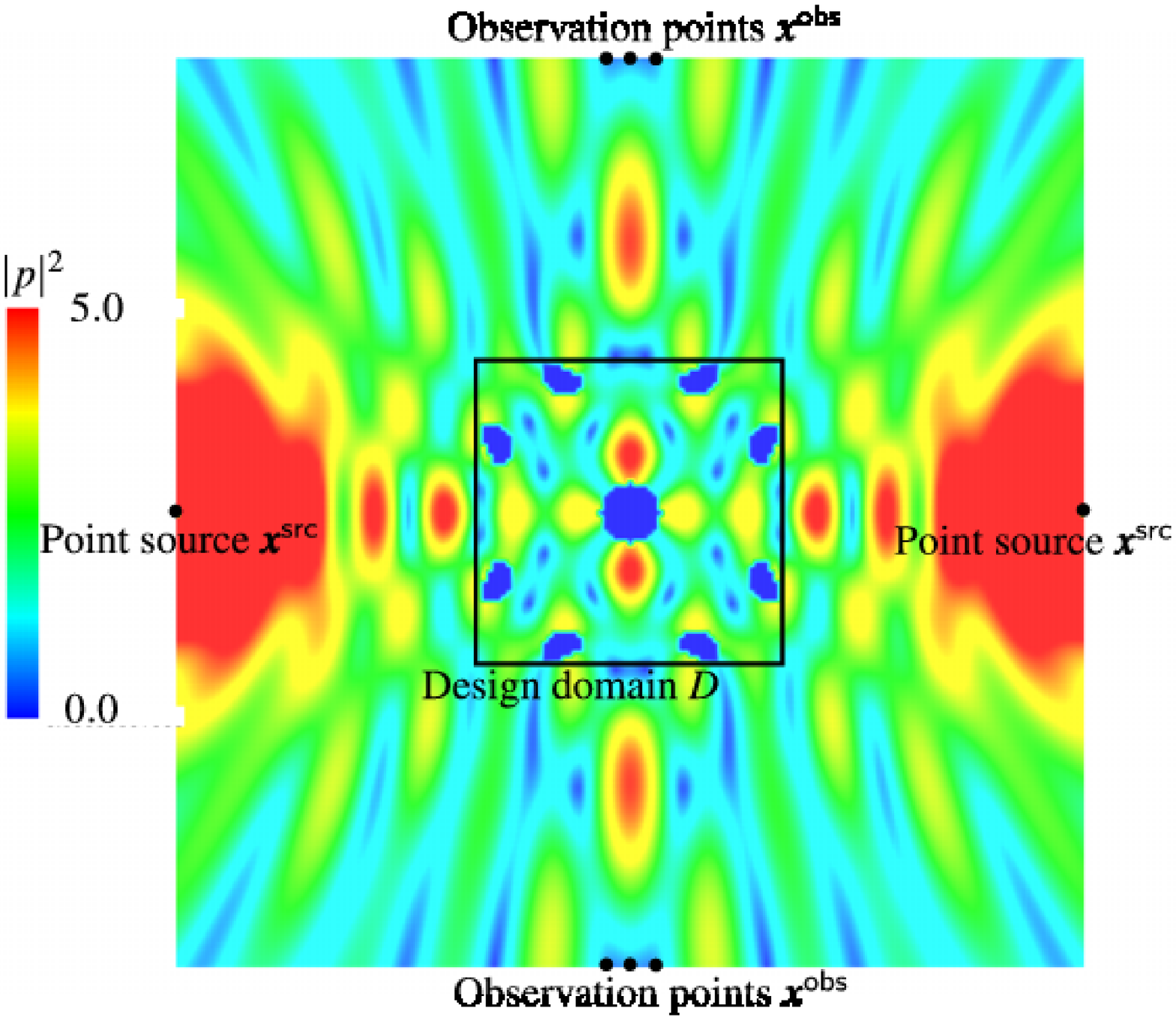}
   \caption{The squared sound norm $|p|^2$ around the design domain $D$ in $x_2=1.25$ when optimal elastic scatterers as Figure~\ref{08tang_mesh} are allocated.}
   \label{10rigidabs}
  \end{center}
  \end{figure}

Figure~\ref{10rigidabs} shows the distribution of the squared sound norm in $x_2=1.25$ plane for the case that optimal elastic scatterers are allocated in the design domain $D$.
One observes that the sound norm on the observation points are reduced.
This can also be confirmed by the objective function, which was $J=42.19$ for the initial configuration in Figure~\ref{075zyouken}, is $J=2.84$ for the optimal configuration.
Thus, the present method can reduce the sound norm on the observation points.

\subsubsection{Soft scatter}
\begin{figure}[h]
  \begin{center}
   \includegraphics[scale=0.12]{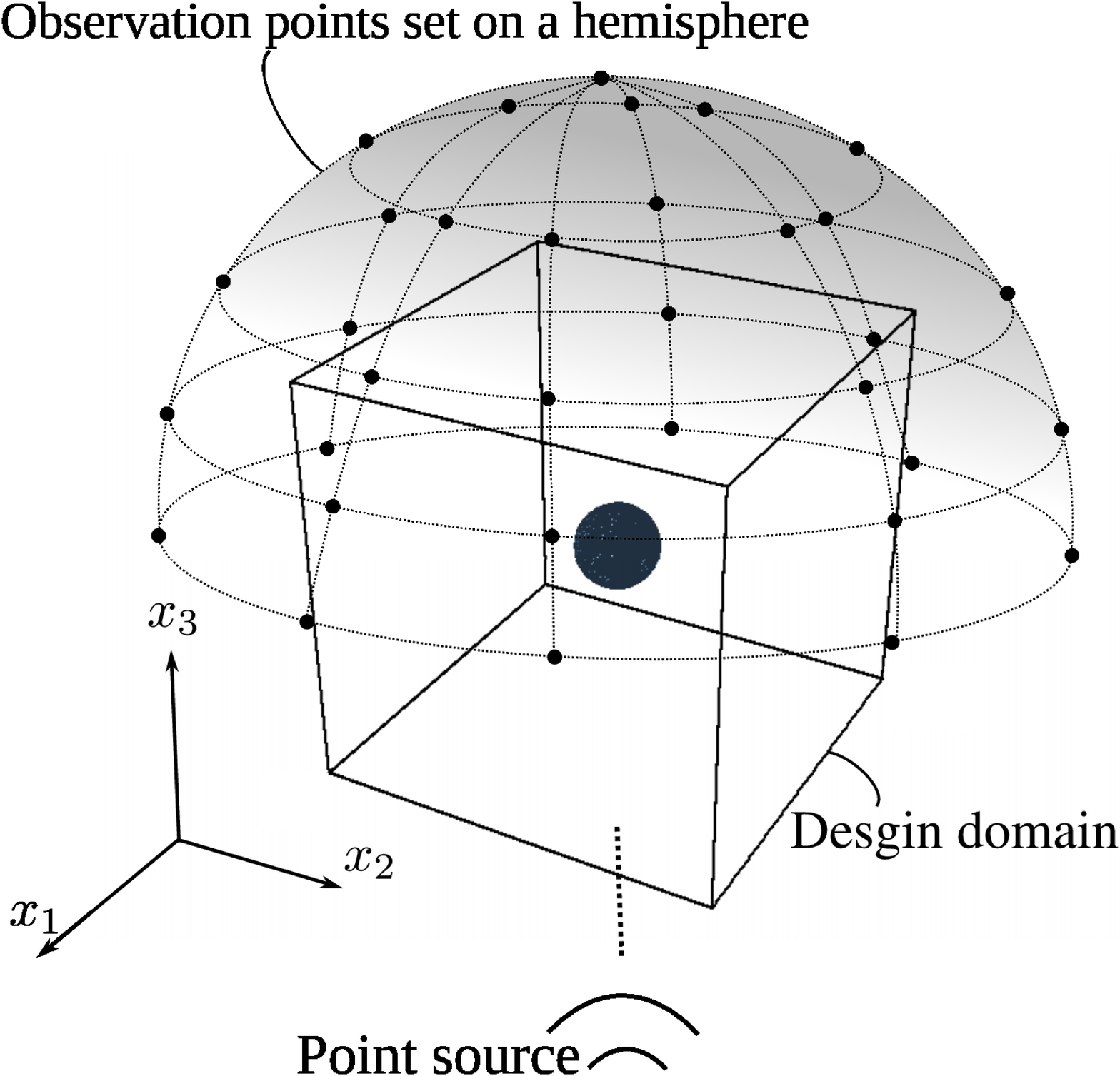}
   \caption{Settings for a topological optimal design for sound scatterer with soft elastic material. In the figure, the initial configuration of the elastic material is also plotted.}
   \label{11zyouken}
  \end{center}
 \end{figure}

We next consider a topology optimisation problem for a soft elastic material to manipulate sound waves.
We use the same design domain and the initial configuration of the elastic scatterer as the one in the previous example, and put a point sound source on $(1.25, 1.25, -2.5)$ whose frequency and intensity are $2\pi$ and $70$, respectively.
The objective function is defined as the sum of the sound norm on 33 observation points set on a hemisphere $\{\boldsymbol{x}~|~\sum_{i=1}^3(x_i-1.25)^2<2.5^2, 1.25\le x_3\}$ (Figure~\ref{11zyouken}).
We here consider a silicone rubber immersed in water.
The material parameters for the silicone rubber is normalised as the density $\rho_\mathsf{s}=0.97$, the Young modulus $E=0.018$ and the Poisson ratio $\nu=0.49$ with the ones for water as the density $\rho_\mathsf{f}=1.0$ and the bulk modulus $\Lambda_\mathsf{f}=1.0$.
Note that the present calculations ignore viscoelastic effect related to silicone rubber.
The viscoelastic effect can be considered by complexifying the Lam\'e constants~\cite{luke1995fluid}, 
which may be addressed in our future publications.
In the optimisation algorithm in Section \ref{lsm}, we divide the design domain $D$ into $100\times 100\times 100$ finite elements, and set $\tau$ in Eqs.~\eqref{timeeq1} and \eqref{timeeq2} as $10^{-4}$ as in the previous example.

\begin{figure}[h]
  \begin{center}
   \includegraphics[scale=0.2]{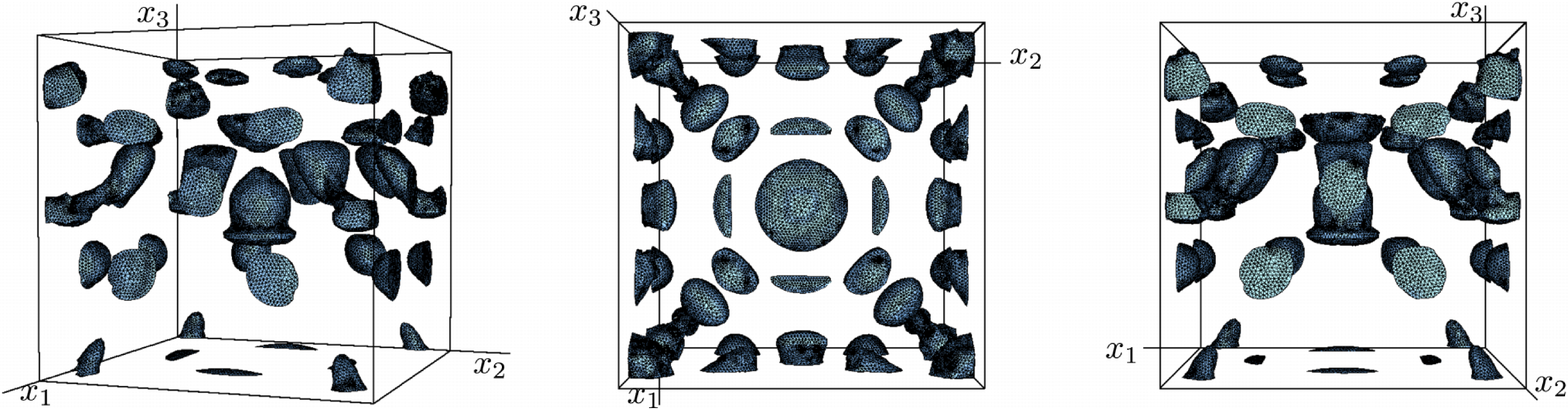}
   \caption{The obtained configuration of silicone rubber.}
   \label{12gomu_mesh}
   \includegraphics[scale=0.2]{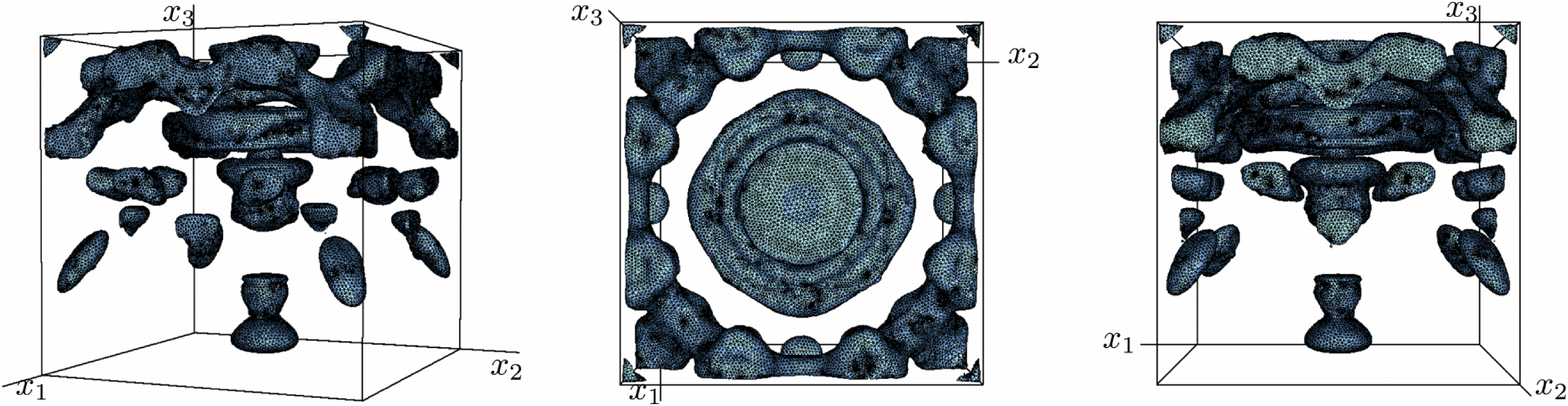}
   \caption{The obtained configuration of rigid material.}
   \label{13rigid_mesh}
  \end{center}
 \end{figure}

\begin{figure}[h]
  \begin{center}
   \includegraphics[scale=0.4]{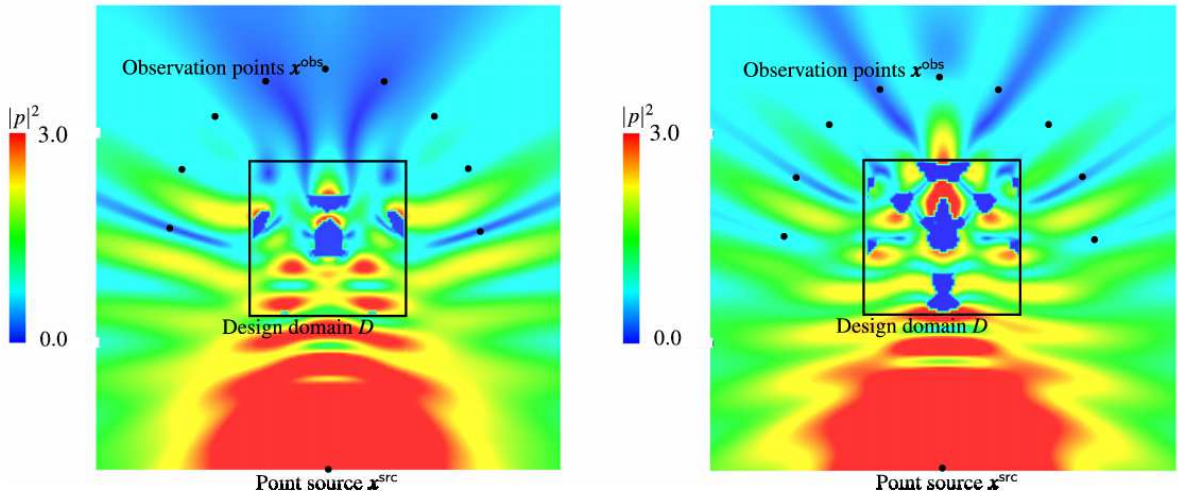}
   \caption{The squared sound norm $|p|^2$ around the design domain $D$ in $x_2=1.25$ when (left:) optimal elastic scatterers (right:) optimal rigid scatterers are allocated.}
   \label{14soft_abs}
  \end{center}
  \end{figure}

We show the obtained configuration in Figure~\ref{12gomu_mesh} 
and sound norm around the optimal scatterers in $x_2=1.25$ in the left figure of Figure~\ref{14soft_abs}.
For comparison, we also show the figures in Figure~\ref{13rigid_mesh} and Figure~\ref{14soft_abs} (right) results of the optimisation problem in which the silicone rubber is replaced by a rigid material.
One observes that the optimal configuration of silicone rubber is different from that of rigid material, 
and, for both cases, the sound norms on the observation points are small.
The objective functions for silicone rubber (resp. rigid material) are reduced from $J=18.49$ (resp. 17.82) to $13.76$ (resp. $5.35$).
These figures shows that, since the material properties of silicone rubber is quite different from those of rigid material, 
such soft materials cannot be designed with conventional topology optimisation with rigid approximation.
  
With these examples, we conclude that the present methodology can efficiently design elastic materials to manipulate sound waves.

\section{Conclusion}
We have developed a new topology optimisation for elastic material to reduce sound level, 
in which a fast BEM-FEM coupled solver is employed to evaluate the topological derivative.
The derivation of the topological derivative for acoustic-elastic coupled problems is described.
We have confirmed that the present topology optimisation method can efficiently design elastic sound scatterers.

In this paper, we have tested the proposed method in pure elastic problem to reduce sound norm defined at observation points.
Applications of the proposed method to viscoelastic design, objective functions defined on boundary such as energy flux on sound absorbing device are, however, still remain to be investigated.
Also, behaviour of the present BEM-FEM solver near resonance (including fictitious resonance) need to be investigated.
In our future publications, we plan to extend the present method to deal with anisotropic and/or the Biot sound absorbers.

\section*{Acknowledgements}
This work was supported by JSPS Grant-in-Aid for Scientific Research (B) (Grant No. 16H04255) and JSPS Grant-in-Aid for Challenging Exploratory Research (Grant No.15K13856).

\bibliographystyle{abbrv}
\bibliography{ref}

\end{document}